\documentclass[12pt,twoside]{amsart}
\usepackage[margin=1in]{geometry}

\usepackage{amsfonts,amsmath,amssymb,amsthm,enumerate,stmaryrd,colonequals,mathrsfs,versions,listings} 
\usepackage[pagebackref=true]{hyperref}
\renewcommand*{\backref}[1]{}
\renewcommand*{\backrefalt}[4]{[{\tiny%
    \ifcase #1 Not cited.%
          \or Cited on page~#2.%
          \else Cited on pages #2.%
    \fi%
    }]}
\usepackage[pdftex]{graphicx}
\usepackage{asymptote}
\usepackage{adjustbox} 

\usepackage[colorinlistoftodos]{todonotes} 

\theoremstyle{plain}

\newcommand{\pf}[1]{\begin{proof}#1 \end{proof}}

\newtheorem{theorem}{Theorem}[section] 
\newcommand{\thm}[1]{\begin{theorem}#1 \end{theorem}}

\newtheorem{theorem-definition}[theorem]{Theorem-Definition} 

\newtheorem{lemma}[theorem]{Lemma} 

\newtheorem{lemma-definition}[theorem]{Lemma-Definition} 

\newtheorem{proposition}[theorem]{Proposition} 
\newcommand{\prop}[1]{\begin{proposition}#1 \end{proposition}}

\newtheorem{proposition-definition}[theorem]{Proposition-Definition} 
\newcommand{\propdefn}[1]{\begin{proposition-definition}#1 \end{proposition-definition}}

\newtheorem{corollary}[theorem]{Corollary} 
\newcommand{\cor}[1]{\begin{corollary}#1 \end{corollary}}

\newtheorem{corollary-definition}[theorem]{Corollary-Definition} 

\newtheorem{observation}[theorem]{Observation} 
\newcommand{\obs}[1]{\begin{observation}#1 \end{observation}}

\newtheorem{exercise}[theorem]{Exercise} 

\newtheorem{problem}[theorem]{Problem} 
\newcommand{\prob}[1]{\begin{problem}#1 \end{problem}}

\newtheorem{question}[theorem]{Question} 
\newcommand{\ques}[1]{\begin{question}#1 \end{question}}

\newtheorem{conjecture}[theorem]{Conjecture} 

\theoremstyle{definition}

\newtheorem{definition}[theorem]{Definition} 
\newcommand{\defn}[1]{\begin{definition}#1 \end{definition}}

\newtheorem{example}[theorem]{Example} 
\newcommand{\ex}[1]{\begin{example}#1 \end{example}}

\theoremstyle{remark}

\newtheorem{remark}[theorem]{Remark} 
\newcommand{\rmk}[1]{\begin{remark}#1 \end{remark}}

\newcommand{\bd}{\mathbf}

\newcommand{\RR}{\mathbb{R}}

\newcommand{\ZZ}{\mathbb{Z}}

\newcommand{\Wlog}{without loss of generality}
\newcommand{\WLOG}{Without loss of generality}

\newcommand{\abs}[1]{\lvert #1 \rvert}
\newcommand{\card}[1]{\lvert #1 \rvert}

\newcommand{\ceil}[1]{\lceil #1 \rceil}

\newcommand{\eps}{\epsilon}

\DeclareFontFamily{U}{wncy}{}
\DeclareFontShape{U}{wncy}{m}{n}{<->wncyr10}{}
\DeclareSymbolFont{mcy}{U}{wncy}{m}{n}
\DeclareMathSymbol{\Sha}{\mathord}{mcy}{"58}

\newcommand{\mcal}{\mathcal}

\newcommand{\texpdf}{\texorpdfstring}

\newcommand{\defeq}{\colonequals} 

\newcommand{\maps}{\colon}

\newcommand{\mathd}[1]{\[#1 \]} 
\newcommand{\matha}[1]{\begin{align*}#1 \end{align*}} 

\newcommand{\belongs}{\subseteq}

\newcommand{\tb}{\textbf}

\newcommand{\set}[1]{\{#1\}} 

\usepackage{xspace}

\newcommand{\Erdos}{Erd\H{o}s\xspace}
\makeatletter
\newcommand{\etale}{\'etal\@ifstar{\'e}{e\xspace}}
\makeatother

\newcommand{\Szemeredi}{Szemer\'{e}di\xspace}

\newcommand{\Holder}{H\"{o}lder\xspace}

\newcommand{\Gyarfas}{Gy\'{a}rf\'{a}s\xspace}

\newcommand{\Szabo}{Szab\'{o}\xspace}

\newcommand{\Dual}{\operatorname{Dual}}
\newcommand{\Color}{\operatorname{Color}}
\newcommand{\Record}{\operatorname{Record}}
\newcommand{\C}{\textup{C}\xspace}
\newcommand{\R}{\textup{R}\xspace}
\newcommand{\G}{\textup{G}\xspace}
\newcommand{\B}{\textup{B}\xspace}
\newcommand{\K}{\textup{K}\xspace}
\newcommand{\X}{\textup{X}\xspace}
\newcommand{\Y}{\textup{Y}\xspace}

\newcommand{\RB}{\textup{RB}\xspace}
\newcommand{\GB}{\textup{GB}\xspace}
\newcommand{\RK}{\textup{RK}\xspace}
\newcommand{\GK}{\textup{GK}\xspace}
\newcommand{\BK}{\textup{BK}\xspace}
\newcommand{\XK}{\textup{XK}\xspace}
\newcommand{\RGB}{\textup{RGB}\xspace}
\newcommand{\RGK}{\textup{RGK}\xspace}
\newcommand{\RBK}{\textup{RBK}\xspace}
\newcommand{\GBK}{\textup{GBK}\xspace}

\newcommand{\RGBK}{\textup{RGBK}\xspace}

\begin{document}



\excludeversion{arxiv}
\includeversion{submit}
\excludeversion{disclaimer}
\excludeversion{markednewpar}




\title{1-color-avoiding paths, special tournaments, and incidence geometry}

\date{\today}

\author{Jonathan Tidor}
\address{Department of Mathematics, Massachusetts Institute of Technology, \mbox{Cambridge, MA 02139}}
\email{jtidor@mit.edu}

\author{Victor Y. Wang}
\address{Department of Mathematics, Massachusetts Institute of Technology, \mbox{Cambridge, MA 02139}}
\email{vywang@mit.edu}

\author{Ben Yang}
\address{Department of Mathematics, Massachusetts Institute of Technology, \mbox{Cambridge, MA 02139}}
\email{beny@math.mit.edu}

\thanks{This version (arXiv v2) acknowledges recent earlier work of Wagner in the case of rainbow-triangle free graphs, which we were unaware of at the time of posting of v1.
In particular, his use of the Gallai--Hasse--Roy--Vitaver theorem allows some technical improvements to one of our results from v1.
However, on a first read, the reader may find it easier to read the self-contained v1.}


\begin{abstract}
We discuss two approaches to a recent question of Loh: must a 3-colored transitive tournament on $N$ vertices have a 1-color-\emph{avoiding} path of vertex-length at least $N^{2/3}$?
This question generalizes the Erd\H{o}s--Szekeres theorem on monotone subsequences.

First, we define three canonical transformations on these tournaments called Color, Record, and Dual.
We use these to establish a reduction to special tournaments with natural geometric and combinatorial properties.
In many cases (including all known tight examples), these tournaments have recursive Gallai decompositions.
Not all relevant tournaments have Gallai decompositions, but those that do satisfy the desired $N^{2/3}$ bound by recent work of Wagner, roughly analogous to earlier work of Fox, Grinshpun, and Pach on a similar \emph{undirected} problem.


Second, we consider the related geometric problem of bounding \emph{slice-increasing} sets $S\subseteq [n]^3$, which---under an additional ordering hypothesis on $S$---was shown by Loh to be equivalent to the original question.
In particular, we establish a rigorous connection from a problem of Szab\'{o} and Tardos, raise a stronger $L^2$-question on slice-counts, and mention a surprising overlap with the joints problem.
\end{abstract}

\maketitle

\tableofcontents

\section{Introduction}

\subsection{Original problems}

We discuss the following Ramsey-type question for 1-color-\emph{avoiding} paths in 3-colored transitive tournaments: Problem 1.1 from a recent paper of Loh \cite{Loh}.
It is perhaps the simplest open generalization of the classical \Erdos--Szekeres theorem on monotone subsequences \cite{ErdosSzekeres}.

\ques{[$L^\infty$-Ramsey]
\label{QUES:Ramsey-for-general-tournaments}

Must every 3-coloring of the edges of the $N$-vertex transitive tournament contain a 1-color-avoiding directed path with at least $N^{2/3}$ vertices?
}

In \cite{Loh}, Loh indirectly used the triangle removal lemma bound to show that there always exists a path of length $\Omega(N^{1/2}e^{\log^* N})$ (beating the trivial bound of $N^{1/2}$).
Later, Wagner \cite{Wagner} obtained the full $N^{2/3}$ bound for \emph{rainbow-triangle free} tournaments by using the recursive Gallai decomposition of such tournaments (cf. Defintion \ref{DEF:undirected-Gallai-RGBK} below), the Gallai--Hasse--Roy--Vitaver theorem, and a suitable analog of the weighted Ramsey's theorem used by Fox, Grinshpun, and Pach \cite{K-free-undirected} for a similar \emph{undirected} problem.

\rmk{
Generally one first asks for the asymptotics, but the exact $N^{2/3}$ threshold is natural for the following reason.
First, the best constructions known achieve exactly $N^{2/3}$.
Furthermore, suppose there were a counterexample $\mcal{T}$ on $N$ vertices, with longest path of length exactly $N^\alpha$, for some real $\alpha < 2/3$.
Then taking the $k$-fold lexicographic ``power'' of $\mcal{T}$ would yield a tournament on $M = N^k$ vertices with longest path of length exactly $(N^\alpha)^k = M^\alpha$---giving arbitrarily large examples with \emph{exponent} $\alpha < 2/3$.
}

We review the definition and properties of the lexicographic product of edge-colored graphs in Definition \ref{DEF:lex-product-tournaments} below.
This leads to the following equivalent question.

\ques{[$L^0$-Ramsey]
\label{QUES:L^0-Ramsey-for-RGB-tournaments}
Given a 3-colored $N$-vertex tournament $\mcal{T}$, must the product of the longest 1-color-avoiding paths in each of the three colors always be at least $N^2$?
}

\pf{[Proof of a priori equivalence of $L^0$- and $L^\infty$- Ramsey problems]
By taking lexicographic products of rotations, one sees that if the $L^\infty$ (single-max) bound holds, then the $L^0$ (geometric mean) bound should hold as well.
The converse holds by pigeonhole.
}

\ex{

Equality holds in Question \ref{QUES:L^0-Ramsey-for-RGB-tournaments}, for instance, when $\mcal{T}$ is a transitive path in one of the three colors, or a lexicographic product thereof.
}

Loh also implicitly introduced a related geometric problem about points in space.

\defn{
\label{DEF:slice-increasing-sets}

Say a set of triples $S\belongs \RR^3$ is \emph{slice-increasing} if every pair of triples in the set is \emph{majority-comparable}, i.e. given $\set{(x,y,z),(x',y',z')} \in \binom{S}{2}$, the coordinate-wise difference $(x'-x,y'-y,z'-z)$ has at least two nonzero coordinates of the same sign.
}

\rmk{
This is equivalent to requiring that on each fixed coordinate-slice, the points are strictly increasing in each of the two remaining coordinates.
}

\ques{
\label{QUES:bounding-slice-increasing-sets}

Must a slice-increasing set $S\belongs [n]^3$ satisfy $\card{S} \le n^{3/2}$?\footnote{Throughout this paper we use the convention $[n] = \set{1,2,\dots,n}$.}
}

To be precise, Loh noted that an ordered set (Definition \ref{DEF:ordered-triples}) is also slice-increasing \cite[Lemma 2.2 and Observation 2.2]{Loh}, and that the question on bounding ordered sets is equivalent to Question \ref{QUES:Ramsey-for-general-tournaments} \cite[Lemma 2.1]{Loh}.
So an affirmative answer to Question \ref{QUES:bounding-slice-increasing-sets} would imply the same for Question \ref{QUES:Ramsey-for-general-tournaments}, while a negative answer for Question \ref{QUES:Ramsey-for-general-tournaments} would imply the same for Question \ref{QUES:bounding-slice-increasing-sets}.

\rmk{
\label{RMK:tensor-power-trick}
Again, the best constructions known achieve exactly $n^{3/2}$ points.
By taking lexicographic powers, a counterexample for Question \ref{QUES:bounding-slice-increasing-sets} would yield arbitrarily large counterexamples with exactly $n^\alpha$ points for some \emph{exponent} $\alpha > 3/2$.
}

We review the lexicographic product of triples in Definition \ref{DEF:lex-product-triples} below.
By considering lexicographic products with rotations, one obtains an \emph{equivalent} problem on unequally-sized grids.
In Section \ref{SEC:geometry-approaches} we will discuss both of these viewpoints.

\ques{
\label{QUES:bounding-unequal-slice-increasing-sets}

Must a slice-increasing set $S\belongs [n_1]\times[n_2]\times[n_3]$ satisfy $\card{S} \le (n_1n_2n_3)^{1/2}$?
}

\rmk{
Given a slice-increasing set, one can extend the $\Color$ map of Definition \ref{DEF:triples-to-tournament} to get a \emph{non-transitive} tournament, and then ask about the longest 1-color-avoiding paths.
However, in this paper we will mainly focus on transitive tournaments.
}

\ex{
\label{EX:standard-slice-increasing-tight-example-sumsets}

Equality holds in Question \ref{QUES:bounding-unequal-slice-increasing-sets}, for instance, when $S$ is the image of a map $A\times B\times C \to \RR^3$ given by $(a,b,c) \mapsto (f(a,b), g(a,c), h(b,c))$, where $f,g,h$ are real-valued, \emph{coordinate-wise strictly increasing}, and \emph{injective} functions defined on $A\times B, A\times C, B\times C$, respectively, and $(n_1,n_2,n_3) = (\card{A}\card{B},\card{A}\card{C},\card{B}\card{C})$.
(Technically, $S$ has $\card{A}\card{B}$ distinct $x$-coordinates, etc. so it is coordinate-wise order-isomorphic to a subset of $[n_1]\times[n_2]\times[n_3]$.)
When $A=B=C = [n^{1/2}]$ with $f=g=h$ given by the linear map $(r,s)\mapsto n^{1/2}r + s$, one recovers a standard lexicographic construction for Question \ref{QUES:bounding-slice-increasing-sets}.
}

\subsection{Structure of paper}

We now summarize our results and the organization of the paper.
In the remainder of the introduction, we build up to Theorem \ref{THM:main-equivalences} reducing the original Question \ref{QUES:L^0-Ramsey-for-RGB-tournaments}---or more precisely, the equivalent Question \ref{QUES:Ramsey-for-RGBK-tournaments}---to Question \ref{QUES:Ramsey-for-canonical-RGBK-tournaments} on special tournaments with a mixture of natural geometric and combinatorial properties (see Definitions \ref{DEF:geometric-tournament} and \ref{DEF:canonical-fixed-points}).
To preserve the flow of the reduction argument, we will usually leave the more routine and technical ingredients to Section \ref{SEC:transformations} and the appendices.

\rmk{
As we will gradually clarify, our reduction process can be viewed as a \emph{saturation} process that \emph{preserves} (global) maximum lengths while \emph{locally} increasing lengths.
Intuitively, this is why we get more structure at the end.
}

Theorem \ref{THM:main-theorem} gives a positive answer to these questions for so-called \emph{undirected-Gallai} tournaments (Definition \ref{DEF:undirected-Gallai-RGBK}).
In v1 of this paper, we were only able to prove the result for \emph{directed-Gallai} tournaments (Definition \ref{DEF:directed-Gallai-RGBK}).
The technical improvement here essentially stems from Wagner's use of the Gallai--Hasse--Roy--Vitaver theorem in \cite{Wagner}.



\rmk{
Wagner \cite{Wagner} actually shows that given any $r$-colored, $N$-vertex \emph{rainbow-triangle free} tournament (not necessarily transitive), there exists---for any fixed integer $0\le s\le r$---a path on $\ceil{N^{s/r}}$ vertices using at most $s$ distinct colors.
We will not say much more about these extensions, because serious difficulties remain in the $(r,s) = (3,2)$ problem (even in the transitive case, which may be easier in view of the history for the $s=1$ problem).
We do note, though, that it does not seem easy to embed the \emph{full} $(r,s) = (3,2)$ problem in these extensions, even after the reductions below.
}

\rmk{
The use of Gallai decompositions, together with Cauchy--Schwarz (or \Holder), essentially forms a \emph{local-to-global} argument.
Indeed, the proof can be interpreted as showing that any counterexample for Question \ref{QUES:Ramsey-for-general-tournaments} (in the case of rainbow-triangle free tournaments, say) would contain a \emph{smaller, more extreme} counterexample.
Perhaps a more robust \emph{local-to-global} argument (not depending on such restrictive Gallai partitioning) could address all tournaments.
}

In Section \ref{SEC:weighted-E-S-and-etc.} we briefly review the \emph{weighting} idea implicit in Wagner's work, which we had independently used in v1 of the present paper.
In particular, for the reader's convenience, we explicitly state a weighted version of \Erdos--Szekeres (Theorem \ref{THM:weighted-Erdos-Szekeres}).
We also mention an application to a problem of \Erdos documented by Steele (Corollary \ref{COR:application-to-Erdos-problem-Steele}), and connections to other problems and interpretations.

In Example \ref{EX:canonical-difficulties} we briefly explain some serious difficulties we have had trying to further extend the proof of Theorem \ref{THM:main-theorem}, even after the reduction in Theorem \ref{THM:main-equivalences}.
Nonetheless, we have included some structural results and Python code on the special tournaments (Definition \ref{DEF:canonical-fixed-points}) in Sections \ref{SEC:Python-code}, \ref{SEC:structure-of-Color-Record-fixed-points}, and \ref{SEC:structure-of-canonical-tournaments}, in case it helps future researchers.

Finally, from a more geometric perspective, we discuss several approaches to the slice-increasing problem (the equivalent Questions \ref{QUES:bounding-slice-increasing-sets} and \ref{QUES:bounding-unequal-slice-increasing-sets}) in Section \ref{SEC:geometry-approaches}, leaving finer details to Section \ref{SEC:slice-increasing-sets} at times in order to maintain coherence.

\subsection{Basic definitions}

As in Loh's paper \cite{Loh}, we will use the classical idea of recording lengths at each vertex.
However, we will also play this $\Record$ map off against two other canonical transformations.
For this, it will help to have the following basic definitions.

\defn{
\label{DEF:ordered-triples}

Call a set of triples $S\belongs \RR^3$ \emph{ordered} if the triples can be listed as $L_1,\dots,L_{\card{S}}$ such that for every $i<j$, the coordinate-wise difference $L_j - L_i$ has at least two strictly positive coordinates.
}

\rmk{
A slice-increasing set $S$ is ordered if and only if the well-defined ``majority-comparable tournament'' on $S$ is acyclic.
Corollary \ref{COR:surfaces-formulation-of-total-order} offers another interpretation.
}

\defn{
Fix an ordered \emph{sequence} of triples $S$.
For any coordinate $c\in \set{x,y,z}$, let $\ell_c$ denote the length of the longest $c$-increasing sub\emph{sequence} of $S$.
}

\ques{
\label{QUES:Ramsey-for-ordered-triples}

Must an ordered set $S\belongs \RR^3$ satisfy $\ell_x \ell_y \ell_z \ge \card{S}^2$?
}

\ques{
\label{QUES:original-triples-problem}

Must an ordered set $S\belongs [n_1]\times [n_2]\times [n_3]$ satisfy $\card{S} \le (n_1n_2n_3)^{1/2}$?
}

\defn{
\label{DEF:RGBK-tournaments}

An \emph{\RGBK-tournament} is a transitive tournament $\mcal{T}$ on ordered vertices $v_1,\dots,v_N$ with each edge $v_i \to v_j$ colored one of the four colors \R, \G, \B, \K.
}

\defn{
Fix an \RGBK-tournament $\mcal{T}$.
For any color class \C, let $\ell_{\mcal{T}}(\C)$ denote the length of the longest \C-colored directed path in $\mcal{T}$. This is abbreviated as $\ell(\C)$ when $\mcal{T}$ is clear.
}

\ques{
\label{QUES:Ramsey-for-RGBK-tournaments}

Must an \RGBK-tournament $\mcal{T}$ satisfy $\ell(\RGK)\cdot \ell(\RBK)\cdot \ell(\GBK) \ge \card{V(\mcal{T})}^2$?
}

Of course, changing a \K-edge to \R, \G, or \B can never increase $\ell(\RGK)$, $\ell(\RBK)$, or $\ell(\GBK)$.
However, when translating between the triples (geometric) and tournament (combinatorial) formulations of the problem, it seems most natural to include a fourth ``wild color'' \K with geometric significance.
In the same vein one may think of \R, \G, \B as ``primary colors'', as we will gradually clarify below.

\subsection{Canonical transformations}

We now define canonical transformations $\Record$, $\Color$, and $\Dual$ as follows.
We also state their basic properties, but leave the simple proofs to Section \ref{SEC:transformations} to preserve the flow of the introduction.

\defn{
\label{DEF:tournament-to-triples}
For an \RGBK-tournament $\mcal{T}$, define $\Record(\mcal{T})\belongs \ZZ_{>0}^3$ to be the following ordered set of $\card{V(\mcal{T})}$ triples.
Each vertex $v$ of $\mcal{T}$ records the triple $(RGK,RBK,GBK)$ where $XYZ$ denotes the length of the longest XYZ-colored path \emph{ending} at $v$.\footnote{These triples are indeed ordered: if $v_i \to v_j$ is colored X, at least two of the sets \RGK, \RBK, \GBK contain X, so $L_i < L_j$ strictly increases in at least two coordinates.}
}

\prop{
\label{PROP:reduction-by-Record}
Transformation $\Record$ sends an \RGBK-tournament $\mcal{T}$ to an ordered set of triples $S\belongs \ZZ_{>0}^3$ with $\card{S} = \card{V(\mcal{T})}$ and $(\ell_x,\ell_y,\ell_z) = (\ell(\RGK),\ell(\RBK),\ell(\GBK))$.
}

\defn{
\label{DEF:triples-to-tournament}

For an ordered set $S\belongs\RR^3$, define $\Color(S)$ to be the following \RGBK-tournament on $\card{S}$ vertices.
For every $i<j$, edge $v_i\to v_j$ is assigned a color via $\R = (+,+,\leqslant\!0)$, $\G = (+,\leqslant\!0,+)$, $\B = (\leqslant\!0,+,+)$, $\K = (+,+,+)$ where a $+$ indicates a \emph{strict} coordinate increase, while a $\leqslant\!0$ indicates a \emph{weak} coordinate decrease.\footnote{By definition of an ordered sequence of triples, the \RGBK color assignment is well-defined.}
}

\rmk{
Suppose one perturbs $S$ to $S'$ via the map $(x,y,z) \mapsto (x-\eps z, y - \eps x, z - \eps y)$.
For small $\eps>0$, the set $S'$ lies in coordinate-general position, and $\Color(S) = \Color(S')$.
}

\prop{
\label{PROP:reduction-by-Color}

Transformation $\Color$ sends an ordered set $S\belongs \RR^3$ to an \RGBK-tournament $\mcal{T}$ with $\card{V(\mcal{T})} = \card{S}$ and $(\ell(\RGK),\ell(\RBK),\ell(\GBK)) = (\ell_x,\ell_y,\ell_z)$.
}

\defn{
\label{DEF:Dual}
Given a tournament $\mcal{T}$, define $\Dual(\mcal{T})$ to be the tournament with all edge directions flipped, but edge colors preserved.
}

\subsection{Reduction tools}

Theorem \ref{THM:main-equivalences} relies on the following basic observations.

\obs{[Reduction by $\Record$]

Let $\mcal{T}$ be an \RGBK-tournament.
If the ordered set $S = \Record(\mcal{T})$ satisfies Question \ref{QUES:Ramsey-for-ordered-triples}, then tournament $\mcal{T}$ satisfies Question \ref{QUES:Ramsey-for-RGBK-tournaments}.
}

\pf{
Proposition \ref{PROP:reduction-by-Record} implies $(\ell_x,\ell_y,\ell_z) = (\ell(\RGK),\ell(\RBK),\ell(\GBK))$, so the bound $\ell(\RGK)\cdot \ell(\RBK)\cdot \ell(\GBK) \ge \card{V(\mcal{T})}^2$ is equivalent to $\ell_x \ell_y \ell_z \ge \card{S}^2$.
}

\obs{[Reduction by $\Color$]

Let $S\belongs \ZZ_{>0}^3$ be an ordered set.
If the \RGBK-tournament $\mcal{T} = \Color(S)$ satisfies Question \ref{QUES:Ramsey-for-RGBK-tournaments}, then set $S$ satisfies Question \ref{QUES:Ramsey-for-ordered-triples}.
}

\pf{
Proposition \ref{PROP:reduction-by-Color} implies $(\ell(\RGK),\ell(\RBK),\ell(\GBK)) = (\ell_x,\ell_y,\ell_z)$, so the bound $\ell_x \ell_y \ell_z \ge \card{S}^2$ is equivalent to $\ell(\RGK)\cdot \ell(\RBK)\cdot \ell(\GBK) \ge \card{V(\mcal{T})}^2$.
}

\obs{[Reduction by $\Dual$]

Let $\mcal{T}$ be an \RGBK-tournament.
If $\Dual(\mcal{T})$ satisfies Question \ref{QUES:Ramsey-for-RGBK-tournaments}, then so does $\mcal{T}$.
}

\pf{
Paths are merely reversed under $\Dual$, so their lengths are preserved.
}

For the statement of Theorem \ref{THM:main-theorem} it will help to have the following definition.

\defn{
\label{DEF:canonically-has-property-P}
Let $\mcal{T}$ be an \RGBK-tournament.
Say $\mcal{T}$ \emph{canonically-almost has} property $\mcal{P}$ if one can apply some finite composition of the maps $\Color\circ\Record$ and $\Dual$ to $\mcal{T}$ to obtain a tournament $\mcal{T}'$ with property $\mcal{P}$.
}

The reduction tools above tell us that an \RGBK-tournament $\mcal{T}$ satisfies Question \ref{QUES:Ramsey-for-RGBK-tournaments} if and only if $\mcal{T}$ \emph{canonically-almost} satisfies Question \ref{QUES:Ramsey-for-RGBK-tournaments}.
Furthermore, one has the following helpful stabilization result, which motivates Definition \ref{DEF:canonical-fixed-points} (canonical tournaments) below.

\prop{
\label{PROP:stabilization-to-canonical-tournaments}
An \RGBK-tournament stabilizes after finitely many applications of $\Color\circ\Record$ and $\Dual\circ\Color\circ\Record\circ\Dual$.
}

\pf{
By Proposition \ref{PROP:Color-Record-RGK-stability}, each transformation weakly increases the number of \K-edges, with equality if and only if the the tournament is fixed under the transformation.
So the process stabilizes after repeatedly applying the two transformations at most as many times as the number of edges of $\mcal{T}$.
}

\subsection{Resulting structure: special tournaments}

We single out the following two classes of tournaments obtained through various compositions of $\Color$, $\Record$, and $\Dual$.

\defn{
\label{DEF:geometric-tournament}

Call an \RGBK-tournament \emph{geometric} if it lies in the image of $\Color$.
}

By Proposition \ref{PROP:image-of-Color-is-blah-transitive}, a tournament is geometric if and only if is transitive in each of the color classes \R, \G, \B, \RGK, \RBK, \GBK (and intersections thereof).
Geometric tournaments are also \emph{interval-connected} (Proposition-Definition \ref{PROPDEF:RGBK-interval-connected}) in each primary color.

\ques{
\label{QUES:Ramsey-for-geometric-RGBK-tournaments}

Must a \emph{geometric} \RGBK-tournament $\mcal{T}$ satisfy $\ell(\RGK)\cdot \ell(\RBK)\cdot \ell(\GBK) \ge \card{V(\mcal{T})}^2$?
}

\defn{
\label{DEF:canonical-fixed-points}

Call an \RGBK-tournament \emph{canonical} if it is fixed under both $\Color\circ\Record$ and $\Dual\circ\Color\circ\Record\circ\Dual$.
}

Canonical tournaments are certainly geometric.
See Appendix \ref{SEC:structure-of-canonical-tournaments} for a full classification of canonical tournaments (Theorem \ref{THM:classification-of-canonical-tournaments}) and other structural properties.
Corollaries \ref{COR:every-vertex-belongs-to-max-paths} and \ref{COR:K-saturation-of-canonical} give some justification for the label ``canonical''.

\ques{
\label{QUES:Ramsey-for-canonical-RGBK-tournaments}

Must a \emph{canonical} \RGBK-tournament $\mcal{T}$ satisfy $\ell(\RGK)\cdot \ell(\RBK)\cdot \ell(\GBK) \ge \card{V(\mcal{T})}^2$?
}

\subsection{Recursive Gallai decomposition}

We start with a classical definition arising from Gallai's classification of rainbow-triangle free \emph{undirected} complete graphs \cite{Gallai}.
We phrase it in a slightly nonstandard way amenable to the ``\K-blind generalizations'' below.

\defn{
\label{DEF:Gallai-in-avoiding-formulation}
Define the class of \emph{Gallai} \RGB-colored undirected complete graphs recursively as follows.
\begin{enumerate}
\item The single-vertex graph is Gallai.

\item Suppose $G$ has a \emph{base decomposition}, meaning a vertex-partition into $m\ge2$ \emph{strictly smaller nonempty} Gallai graphs $H_1,\dots,H_m$, where the edges between two distinct blocks $H_i,H_j$ use \emph{at most one} of the colors \R, \G, \B, and the edges between the various blocks $H_1,\dots,H_m$ \emph{in total} use \emph{at most two} of the colors \R, \G, \B.
Then $G$ is Gallai.
\end{enumerate}
Furthermore, given a Gallai graph, call any such recursive sequence of vertex-partitions a \emph{(recursive) Gallai decomposition}.
}

We now explicitly recall the original motivation for the previous definition.

\thm{[Gallai \cite{Gallai}]
\label{THM:Gallai-original-RGB-colored}

Let $G$ be a rainbow-triangle free \RGB-colored undirected complete graph.
Then $G$ is Gallai.
In fact, $G$ is always disconnected in one of \R, \G, \B, and if $G$ is disconnected in color \X, then the \X-connected components form a base decomposition of $G$.
}

See \Gyarfas--Simonyi \cite[Theorem A]{UndirectedRainbowFreeCharacterization} for a clearer isolation of the result.
We now give four natural related notions, which define the context of Theorem \ref{THM:main-theorem}.

\defn{
\label{DEF:undirected-Gallai-RGBK}
Call an \RGBK-tournament \emph{undirected-Gallai} if it has an \emph{undirected, \K-blind} Gallai decomposition.
This is defined recursively, word-for-word as in Definition \ref{DEF:Gallai-in-avoiding-formulation}, without additional conditions on \K-edges.\footnote{Essentially, \K-edges can ``go anywhere''.}
}

\propdefn{
\label{PROPDEF:undirected-gallai-vs-rainbow-triangle-free}

An \RGBK-tournament $\mcal{T}$ is undirected-Gallai if and only if it is \emph{morally rainbow-triangle free}, meaning that one can assign each \K-edge of $\mcal{T}$ a new primary color among \RGB to get a rainbow-triangle free tournament $\mcal{T}'$.
}

\pf{
A morally rainbow-triangle free tournament $\mcal{T}$ is automatically undirected-Gallai by Theorem \ref{THM:Gallai-original-RGB-colored} and the \K-blindness of Definition \ref{DEF:undirected-Gallai-RGBK}.
Conversely, one proves recursively that an undirected-Gallai tournament is morally rainbow-triangle free: note that at each level, the base decomposition is, modulo \K-edges, the blowup of a 2-colored graph.
}

Sometimes it is more natural to require direction, as follows.
But the reader can already safely skip to Theorems \ref{THM:main-equivalences} and \ref{THM:main-theorem} below for our concrete results.

\defn{
\label{DEF:directed-Gallai-RGBK}
Call an \RGBK-tournament \emph{directed-Gallai} if it has a \emph{directed, \K-blind} Gallai decomposition, meaning an undirected, \K-blind decomposition where the base decomposition $G = H_1\sqcup \dots\sqcup H_m$ at each step must be \emph{directed}, in the sense that for $i<j$, the vertices in $H_i$ are all directed towards the vertices in $H_j$.
}

\ex{[\K-blindness]
A transitive rainbow-triangle with edges \R, \G, \K is directed-Gallai (and undirected-Gallai), while a transitive rainbow-triangle with edges \R, \G, \B is not directed-Gallai (nor undirected-Gallai).
}

\defn{
\label{DEF:morally-K-free-RGBK}
Call an \RGBK-tournament $\mcal{T}$ \emph{morally \K-free} if one can assign each \K-edge of $\mcal{T}$ a new \emph{primary color} among \RGB to get a \K-free \emph{geometric} tournament $\mcal{T}'$.
}

By Proposition \ref{PROP:image-of-Color-is-blah-transitive}, a \K-free tournament is geometric if and only if it is transitive in all color combinations,  or equivalently if it is single-color-transitive and rainbow-triangle free.
For such tournaments we may apply Theorem \ref{THM:Gallai-original-RGB-colored} as follows.

\prop{
\label{PROP:geometric-K-free-is-directed-Gallai}

Any \K-free geometric \RGBK-tournament is directed-Gallai.
}

\pf{
Since the tournament $\mcal{T}$ is \RGB-colored and rainbow-triangle free, it has an a priori \emph{undirected, \K-blind} decomposition by Theorem \ref{THM:Gallai-original-RGB-colored}.
However, suppose for every base decomposition in the recursion we use connected components in a disconnected color, as allowed by the second clause of Theorem \ref{THM:Gallai-original-RGB-colored}.
Then we in fact obtain a \emph{directed, \K-blind} decomposition by the interval-connectivity observation of Proposition-Definition \ref{PROPDEF:RGBK-interval-connected}.
}

We now relate the previous three definitions.

\prop{
\label{PROP:morally-K-free-implies-directed-Gallai}
If an \RGBK-tournament is directed-Gallai, it is undirected-Gallai.
If it is morally \K-free, it is directed-Gallai and undirected-Gallai.
}

\pf{[Proof that directed-Gallai implies undirected-Gallai]

A directed decomposition is automatically an undirected decomposition as well.
}


\pf{[Proof that morally \K-free implies directed-Gallai]

By Definition \ref{DEF:morally-K-free-RGBK}, take a \K-free geometric tournament $\mcal{T}'$ obtained from a morally \K-free tournament $\mcal{T}$ by changing each \K-edge to a primary color.
By Proposition \ref{PROP:geometric-K-free-is-directed-Gallai}, $\mcal{T}'$ has a directed, \K-blind Gallai decomposition.
By \K-blindness it is in fact a directed, \K-blind Gallai decomposition for $\mcal{T}$ as well.
}




\subsection{Statement of main tournament results}


\thm{
\label{THM:main-equivalences}

Questions \ref{QUES:Ramsey-for-general-tournaments}, \ref{QUES:L^0-Ramsey-for-RGB-tournaments}, \ref{QUES:Ramsey-for-ordered-triples}, \ref{QUES:original-triples-problem}, \ref{QUES:Ramsey-for-RGBK-tournaments}, \ref{QUES:Ramsey-for-geometric-RGBK-tournaments} and \ref{QUES:Ramsey-for-canonical-RGBK-tournaments} are all equivalent.
}

\pf{
Reduction tools $\Color$ and $\Record$ show that Question \ref{QUES:Ramsey-for-ordered-triples} (Ramsey for triples) is equivalent to Question \ref{QUES:Ramsey-for-RGBK-tournaments} (Ramsey for tournaments).
An affirmative answer to Question \ref{QUES:Ramsey-for-ordered-triples} would directly imply the same for Question \ref{QUES:original-triples-problem} (bounding triples).
Conversely, an affirmative answer to the latter would imply the same for the former, by using $\Record\circ\Color$, which sends an ordered set $S\belongs \RR^3$ to an ``efficiently-packed'' ordered set $S'\belongs [\ell_x]\times[\ell_y]\times[\ell_z]$.

It remains to show equivalence of the tournament questions.
We already showed Questions \ref{QUES:Ramsey-for-general-tournaments} and \ref{QUES:L^0-Ramsey-for-RGB-tournaments} are equivalent.
Question \ref{QUES:L^0-Ramsey-for-RGB-tournaments} is equivalent to the superficially more general Question \ref{QUES:Ramsey-for-RGBK-tournaments}, since changing \K's to primary colors cannot increase 1-color-avoiding path lengths, locally or globally.
Finally, Question \ref{QUES:Ramsey-for-canonical-RGBK-tournaments} is equivalent to the ostensibly more general Questions \ref{QUES:Ramsey-for-RGBK-tournaments} and \ref{QUES:Ramsey-for-geometric-RGBK-tournaments}.
This immediately follows from Proposition \ref{PROP:stabilization-to-canonical-tournaments}, which says that every \RGBK-tournament is canonically-almost (Definition \ref{DEF:canonically-has-property-P}) canonical.
}

In view of Theorem \ref{THM:main-equivalences} reducing Question \ref{QUES:Ramsey-for-general-tournaments} to Question \ref{QUES:Ramsey-for-canonical-RGBK-tournaments}, our following (main) theorem is in some sense most immediately useful in the case that $\mcal{T}$ is canonical (Definition \ref{DEF:canonical-fixed-points}).
However, we state it in natural generality for convenience of the reader.

\thm{
\label{THM:main-theorem}

Let $\mcal{T}$ be an \RGBK-tournament.
If $\mcal{T}$ is canonically-almost morally \K-free, canonically-almost directed-Gallai, canonically-almost undirected-Gallai, \emph{or} canonically-almost morally rainbow-triangle free, then $\mcal{T}$ satisfies Question \ref{QUES:Ramsey-for-RGBK-tournaments}.
}

\pf{
By Proposition \ref{PROP:morally-K-free-implies-directed-Gallai} and Proposition-Definition \ref{PROPDEF:undirected-gallai-vs-rainbow-triangle-free}, we immediately reduce to the case of rainbow-triangle free  tournaments.
The result now follows by \cite[Theorem 1.6]{Wagner}.
}


\rmk{[Failure of naive reduction to \K-free geometric case]
If one starts with a geometric or even canonical tournament, and replaces each \K with \R (say) to get a merged color \RK, then one has transitivity in all combinations of \RK, \G, \B \emph{except} \GB---so not quite a \K-free geometric tournament in general.
}

\ex{
\label{EX:canonical-difficulties}

Not all canonical tournaments are morally \K-free, directed-Gallai, undirected-Gallai, or morally rainbow-triangle free.
To see this, it suffices by Proposition \ref{PROP:morally-K-free-implies-directed-Gallai} and Proposition-Definition \ref{PROPDEF:undirected-gallai-vs-rainbow-triangle-free} to find canonical tournaments that are not undirected-Gallai.
In Section \ref{SEC:Python-code} we give an 8-vertex example found by random search.
}

\rmk{
It may still be possible to completely reduce Question \ref{QUES:Ramsey-for-general-tournaments} to the situation of Theorem \ref{THM:main-theorem}, using other transformations.
For example, the set of canonical tournaments is \emph{not} stable under lexicographic product, even though the set of geometric tournaments is.
}


\section{Details of transformations on tournaments and triples}
\label{SEC:transformations}

\subsection{Record: from tournaments to triples}

For the definition and statement of basic properties, see Definition \ref{DEF:tournament-to-triples} and Proposition \ref{PROP:reduction-by-Record} above.

\pf{[Proof of Proposition \ref{PROP:reduction-by-Record}]
By definition, $S$ has $\card{V(\mcal{T})}$ triples.
Now consider any $x$-increasing sequence of $\ell$ triples $L_{i_1}<\dots<L_{i_\ell}$ in $S$.
Then $x_{i_\ell} \ge \ell$ since triples in $S$ have positive integer coordinates.
But by definition, $x_{i_\ell}$ simply records the length of the longest \RGK-path ending at the $i_\ell$th vertex of $\mcal{T}$.
In particular, $\ell(\RGK) \ge \ell$, so $\ell(\RGK) \ge \ell_x$.

Conversely, let $j_1\to \dots \to j_\ell$ be an \RGK-path in $\mcal{T}$ of length $\ell$.
Then $x_{j_1} < \dots < x_{j_\ell}$ (because any \RGK-path ending at $j_i$ can be extended via the \RGK-edge $j_i\to j_{i+1}$ for $i < \ell$), so $L_{j_1}<\dots<L_{j_\ell}$ is an $x$-increasing sequence of $\ell$ triples.
Thus $\ell_x \ge \ell(\RGK)$.

We conclude that $\ell_x = \ell(\RGK)$.
Similarly, $\ell_y = \ell(\RBK)$ and $\ell_z = \ell(\GBK)$.
}

\rmk{
The equality is a little subtle.
It is certainly not true that $i_1 \to \dots \to i_\ell$ is an \RGK-path in $\mcal{T}$ whenever $L_{i_1}<\dots<L_{i_\ell}$ is an $x$-increasing sequence of triples in $S$.
For a more striking instance of the subtlety of the $\Record$ map, see Proposition \ref{PROP:Color-Record-RGK-stability} below.
}

\subsection{Color: from triples to tournaments}

For the definition and statement of basic properties, see Definition \ref{DEF:triples-to-tournament} and Proposition \ref{PROP:reduction-by-Color} above.

\pf{[Proof of Proposition \ref{PROP:reduction-by-Color}]
By definition, $\mcal{T}$ has $\card{S}$ vertices.
Furthermore, by definition, an edge $i\to j$ is \RGK-colored if and only if $L_j - L_i$ has strictly positive $x$-coordinate.
So a directed path $i_1 \to \dots \to i_\ell$ in $\mcal{T}$ is \RGK-colored if and only if the sequence of triples $L_{i_1}<\dots<L_{i_\ell}$ in $S$ is $x$-increasing.
Thus $\ell(\RGK) = \ell_x$.
Similarly, $\ell(\RBK) = \ell_y$ and $\ell(\GBK) = \ell_z$.
}

We now classify geometric \RGBK-tournaments (Definition \ref{DEF:geometric-tournament}).
What really matters for Theorem \ref{THM:main-theorem} is the easy (only if) direction, but we have included both directions for conceptual clarity.

\prop{
\label{PROP:image-of-Color-is-blah-transitive}
An \RGBK-tournament is geometric if and only if it is transitive in each of the color classes \R, \G, \B, \RGK, \RBK, \GBK.
}

\pf{
If $\mcal{T}$ is geometric, then it is easy to check each of the transitivity conditions, because an \R-edge, for instance, is equivalent to a weak $z$-decrease, whereas a \GBK-edge is equivalent to a strict $z$-increase.
Conversely, suppose an abstract tournament $\mcal{T}$ satisfies the transitivity conditions.
Then we may construct $x_1,\dots,x_N$ inductively.
Indeed, once $x_1,\dots,x_{k-1}$ have been constructed, transitivity in \RGK and \B allows $x_k$ to either be sandwiched between some two unique neighboring terms of $x_1,\dots,x_{k-1}$, or else uniquely placed at one of the two ends.
This completes the induction.
After analogously creating $y_1,\dots,y_N$ and $z_1,\dots,z_N$, the resulting list of triples $L_i = (x_i,y_i,z_i)$ must be ordered, and it maps to $\mcal{T}$ under $\Color$.
To see this, note that for $i<j$, the difference-type of $L_j - L_i$ matches the color of $i \to j$ in $\mcal{T}$ on each of the coarse comparisons \RGK vs. \B ($x$-coordinate), \RBK vs. \G ($y$-coordinate), and \GBK vs. \R ($z$-coordinate), all by construction.
So they exactly match.
}

\rmk{
In the converse, all that matters is the order-isomorphism classes of (i.e. permutations corresponding to) the sequences $x_1,\dots,x_N$; $y_1,\dots,y_N$; and $z_1,\dots,z_N$.
}

Concretely, $\mcal{T}$ is geometric if and only if it is transitive in each of \R, \G, \B, \K, and every rainbow-triangle has color \K on the ``long'' edge.
One may also think of transitivity in \R and \GBK as follows: if $i < j < k$ with $ij, jk$ colored \R, then certainly $ik$ is also \R.
But a partial converse holds as well: if $ik$ is \R, then one of $ij, jk$ must be \R.
We can use the partial converse as follows.

\propdefn{
\label{PROPDEF:RGBK-interval-connected}

Let $\mcal{T}$ be a geometric \RGBK-tournament.
Then $\mcal{T}$ is \emph{interval-connected} in each color class $\C \in \set{\R, \G, \B, \RGK, \RBK, \GBK}$ (and unions thereof), meaning that every connected component in class \C consists of a full interval of indices $[i,j]$.
}

\pf{
By transitivity in the complement $\RGBK\setminus\C$, we know that for any \C-colored edge $i\to j$, the interval $[i,j]$ is \C-connected: in fact, each $k\in (i,j)$ is directly \C-connected to one of the endpoints $i,j$.
But then for any \emph{undirected} \C-path $i_0,i_1,\dots,i_\ell$ (with $\ell\ge1$), the union of the (undirected) intervals $[i_0,i_1],\dots,[i_{\ell-1},i_\ell]$ must be \C-connected as well.
By the intermediate value theorem, this union of intervals contains the (undirected) interval $[i_0,i_\ell]$.
We conclude that any \C-component is a full interval: one can take $i_0$ to be the smallest vertex of the component, and $i_\ell$ to be the largest.
}

\subsection{Dual of a tournament}

See Definition \ref{DEF:Dual} above.

\subsection{Lexicographic product of tournaments}

\defn{[Cf. {\cite[Definition 1.2]{K-free-undirected}}]
\label{DEF:lex-product-tournaments}

Given \RGBK-tournaments $\mcal{T}_1,\mcal{T}_2$, define the \emph{lexicographic product} $\mcal{T}_1\otimes \mcal{T}_2$ to be the \RGBK-tournament on lexicographically-ordered vertices $(v_1,v_2)\in \mcal{T}_1\times \mcal{T}_2$ such that the edge $(v_1,v_2)\to (w_1,w_2)$ inherits the color of $v_1\to w_1$ unless $v_1 = w_1$, in which case it inherits the color of $v_2\to w_2$.
}

This can also be interpreted in terms of blowups, where one places a copy of $\mcal{T}_2$ at each vertex of $\mcal{T}_1$.
We document the following basic properties of lexicographic products of tournaments.

\prop{
\label{PROP:lex-product-tournament-properties}

The lexicographic product of tournaments is associative but noncommutative.
In particular, powers are uniquely defined.
Also, for tournaments $\mcal{T}_1,\mcal{T}_2$ one has $\card{V(\mcal{T}_1\otimes \mcal{T}_2)} = \card{V(\mcal{T}_1)}\card{V(\mcal{T}_2)}$, and $\ell_{\mcal{T}_1\otimes \mcal{T}_2}(\C) = \ell_{\mcal{T}_1}(\C) \ell_{\mcal{T}_2}(\C)$ for any color class \C.
}

\subsection{Lexicographic product of triples}

\defn{[Cf. {\cite[proof of Theorem 5 for $d=2$, before Lemma 6]{SzaboTardos}}]
\label{DEF:lex-product-triples}

Given arbitrary sets of triples $S_1,S_2$ \emph{endowed with embeddings} $S_1\belongs [\ell_1]\times[m_1]\times[n_1]$ and $S_2\belongs [\ell_2]\times[m_2]\times[n_2]$, define the \emph{lexicographic product} $S_1\otimes S_2$ to be the embedded set $S\belongs [\ell_1\ell_2]\times [m_1m_2]\times [n_1n_2]$ consisting of triples of the form $(\ell_2 (x_1-1) + x_2, m_2 (y_1-1) + y_2, n_2 (z_1-1) + z_2)$, with $(x_1,y_1,z_1)\in S_1$ and $(x_2,y_2,z_2)\in S_2$.
}


\prop{
\label{PROP:lex-product-triple-properties}

The lexicographic product of \emph{embedded} sets of triples is associative but noncommutative.
In particular, powers are uniquely defined.
Also, for embedded sets of triples $S_1,S_2$ one has $\card{S_1\otimes S_2} = \card{S_1}\card{S_2}$.
If $S_1,S_2$ are furthermore slice-increasing (resp. ordered), then $S_1\otimes S_2$ is slice-increasing (resp. ordered) as well.
}

\rmk{
Slice-increasing sets are actually stable under a less canonical generalization of Definition \ref{DEF:lex-product-triples}.
Indeed, given slice-increasing sets $S_1,S_2\belongs \RR^3$, let $\alpha,\beta,\gamma\maps \RR^2 \to \RR$ be coordinate-wise strictly increasing functions.
Suppose that $\alpha,\beta,\gamma$ are \emph{injective} on $X_1\times X_2, Y_1\times Y_2, Z_1\times Z_2$, respectively, where for $i=1,2$, we define $X_i$ (resp. $Y_i$; $Z_i$) to be the coordinate projection of $S_i$ onto the $x$-axis (resp. $y$-axis; $z$-axis).
Consider the image $S$ of the map $S_1\times S_2\to \RR^3$ given by $(x_1,\dots,x_2,\dots) \mapsto (\alpha(x_1,x_2),\dots)$.
Then $S$ is slice-increasing.
}

\section{Weighted \texpdf{\Erdos}{Erdos}--Szekeres and friends}
\label{SEC:weighted-E-S-and-etc.}

\subsection{Statement and proof}

It will be convenient to make the following definition.

\defn{
\label{DEF:RBK-tournament-and-geometric}

An \emph{\RBK-tournament} is a \G-free \RGBK-tournament.
Call an \RBK-tournament \emph{geometric} if it is geometric as an \RGBK-tournament, or equivalently if it is transitive in each of the color classes \R, \B, \RK, \BK.
}

The following \emph{weighting idea} is implicit in the work of Wagner \cite{Wagner} (and was also used in v1 of the present paper), but for the reader's convenience we state and prove it explicitly.

\thm{[Weighted \Erdos--Szekeres and \RBK-tournaments generalization]
\label{THM:weighted-Erdos-Szekeres}

Let $\mcal{H}$ be an \RBK-tournament (resp. geometric \RBK-tournament) on $M$ vertices.
Consider any \emph{nonnegative} reals $B_1,\dots,B_M$ and $R_1,\dots,R_M$.
Let $B \defeq \max_{P} \sum_{i\in P} B_i$ and $R \defeq \max_{Q} \sum_{j\in Q} R_j$, where $P$ and $Q$ run over paths (resp. cliques) in $\mcal{H}$ of color \BK and \RK, respectively.
Then
\mathd{
B\cdot R \ge \sum_{i=1}^{M} B_i\cdot R_i.
}
}

\rmk{
For $B_i = R_j = 1$ we recover the unweighted version.
}


\pf{
First, throw out all vertices $i$ with $B_i\cdot R_i = 0$.
Now by scaling and rational approximation, reduce to the case of distinct \emph{positive integer} weights $B_i,R_j$.
In this case we can give a combinatorial reformulation.
Let $\mcal{G}$ be the \RBK-tournament (resp. geometric \RBK-tournament) on $\sum_{i=1}^{M} B_i\cdot R_i$ vertices formed by blowing up vertex $i\in \mcal{H}$ into the standard lexicographic tournament on $[B_i]\times[R_i]$, explicitly with the following coloring: the directed edge $(u,v) < (u',v')$ is colored \B if $u<u'$, and \R otherwise (if $u=u'$ but $v<v'$).

Then we merely wish to show that the product of the respective lengths $B$ and $R$ of the longest \BK- and \RK- paths of $\mcal{G}$ is at least $\sum_{i=1}^{M} B_i\cdot R_i = \card{V(\mcal{G})}$.
This immediately follows from the usual \Erdos--Szekeres (or more precisely, the \RBK-tournaments generalization, proven as usual using $\Record$ for color classes \BK and \RK).
}

\subsection{Application to a problem of \texpdf{\Erdos}{Erdos}}

At the end of his review \cite[Section 12]{Steele}, Steele mentions a ``question posed by \Erdos (1973) for which there seems to have been no progress'':

\prob{
Given $x_1,\dots,x_n$ distinct real numbers determine $\max_M \sum_{i\in M} x_i$ over all subsets $M\belongs [n]$ of indices $i_1<\dots<i_k$ such that $x_{i_1},\dots,x_{i_k}$ is monotone.
}

\cor{
\label{COR:application-to-Erdos-problem-Steele}

In the above situation, $\max_M \sum_{i\in M} x_i \ge (\sum_i \max(x_i,0)^2)^{1/2}$, if we use the convention that the empty sum is $0$.
}

\pf{
Construct the usual transitive \RB-tournament on vertices $v_1,\dots,v_n$, with $v_i \to v_j$ colored \R (say) if $x_i < x_j$, and \B if $x_i > x_j$.
If $M$ maximizes $\sum_{i\in M} x_i$, then $x_i \ge 0$ for all $i\in M$, so
\mathd{
\left(\max_M \sum_{i\in M} x_i\right)^2
= \left(\max_M \sum_{i\in M} \max(x_i,0)\right)^2
\ge 
\sum_{i=1}^n \max(x_i,0)^2
}
by Theorem \ref{THM:weighted-Erdos-Szekeres}, as desired.
}

\subsection{Tournaments transitive in every color combination}

In this section, we briefly discuss the geometric \K-free case of Theorem \ref{THM:main-theorem}.
It has clean geometric interpretations and connections to other problems.
For instance, the result has a \emph{local-to-global} interpretation: suppose every triangle avoids one of the three colors \RGB.
Then by the geometric transitivity properties, there is a $\ceil{N^{2/3}}$-vertex sub-tournament---not just a path---avoiding one of the three colors \RGB (cf. Ramsey discussion in Section \ref{SUBSEC:Ramsey-and-Szabo-Tardos}).

\subsubsection{``\K-flatness''}
\label{SUBSUBSEC:K-flatness}

Observe the following ``\K-flat'' tight construction for Question \ref{QUES:Ramsey-for-ordered-triples}.

\ex{
Start with the list $\set{(0,0,0),(1,1,-2)}$ and take the lexicographic product with its two cyclic rotations, in any order.
(It is easy to check that the product is \K-free under $\Color$.)
The product may be embedded in the ``\K-flat'' plane $x+y+z=0$.
This can be viewed as a ``\K-flat'' version of the usual equality case for Question \ref{QUES:original-triples-problem}, where one starts with the list $\set{(0,0,0),(1,1,0)}$ instead, and the product is no longer \K-free under $\Color$.
}

\rmk{
Let $S \belongs [n]^3$ be a sequence of triples such that for $i<j$, the triples $(x_i,y_i,z_i)$ and $(x_j,y_j,z_j)$ strictly increase in two coordinates \emph{and weakly decrease in the third coordinate}.
Then Theorem \ref{THM:main-theorem} only directly gives a bound of $\card{S}\leq n^{3/2}$.
However, using the fact that $S$ is constrained to $[n]^3$, one can in fact get a linear bound: see Proposition \ref{PROP:linear-bound-on-K-free-slice-increasing-sets} below.
}




\section{Geometric approaches: triples and grids}
\label{SEC:geometry-approaches}

In this section, we explore various geometric approaches to bounding the size of a constrained slice-increasing set (Questions \ref{QUES:bounding-slice-increasing-sets} and \ref{QUES:bounding-unequal-slice-increasing-sets}).
Although for the original transitive tournament problem one really cares about ordered sets (Question \ref{QUES:original-triples-problem}), most of the ideas below apply more generally to slice-increasing sets.
We will try to clarify when we believe one may be able to truly leverage the ordering, such as in Section \ref{SUBSEC:L^2-slice-counts-approach}.
At a higher level, note that the reduction in Theorem \ref{THM:main-equivalences} only applies to ordered triples and their associated tournaments; for geometric facts related to the reduction, we direct the interested reader to Theorem \ref{THM:classification-of-Color-Record-fixed-points}, Theorem \ref{THM:classification-of-canonical-tournaments}, and Corollary \ref{COR:min-above,max-below}.

\subsection{Outline of section}

In Section \ref{SUBSEC:grids-intro} we introduce two-dimensional ``grid views'' (Definition \ref{DEF:projection-and-grid-view-of-triples}) for visualizing sets $S\belongs \ZZ^3$, to be used throughout Section \ref{SEC:geometry-approaches}.
In Section \ref{SUBSEC:Ramsey-and-Szabo-Tardos} we give a Ramsey perspective to Question \ref{QUES:bounding-slice-increasing-sets}, and show (in Theorem \ref{THM:Szabo-Tardos-transfers-to-slice-increasing-bounds}) that any nontrivial bound on a problem of \Szabo and Tardos would transfer over to the problem of bounding slice-increasing sets.
In Section \ref{SUBSEC:L^2-slice-counts-approach} we suggest a natural $L^2$-question on slice-counts with vaguely possible connections to representation theory; again, any nontrivial bound would transfer over to the problem of bounding slice-increasing sets.
In Section \ref{SUBSEC:joints-problem-inspiration} we observe a surprising overlap between tight examples for Question \ref{QUES:bounding-slice-increasing-sets} and the joints problem \cite[Theorem 1.1]{GuthKatz}.

When it would take us too far astray to discuss a given approach in too much greater detail, we will refer accordingly to Appendix \ref{SEC:slice-increasing-sets} for the interested reader.

\subsection{Grids}
\label{SUBSEC:grids-intro}

The ``grid view'' hinges on \cite[Observation 2.1]{Loh} that no two triples in a slice-increasing set agree in exactly two coordinates.

\defn{[Cf. {\cite[bipartite graph construction immediately after Observation 2.1]{Loh}}]
\label{DEF:projection-and-grid-view-of-triples}

Let $S\belongs \ZZ^3$ be any set of triples for which no two triples agree in \emph{exactly} two coordinates.
Let $S_z\belongs \ZZ^2$ denote the $xy$-plane \emph{projection} of $S$, so that $S\to S_z$ is injective.
To obtain the \emph{$xy$-grid view} of $S$, we leave a square $(x,y)\in \ZZ^2$ \emph{empty} or \emph{unlabeled} if $(x,y)\notin S_z$, and otherwise fill it in with the unique \emph{label} $z\in \ZZ$ such that $(x,y,z) \in S$.
}

\rmk{
One can do the same for the $xz$-plane and $yz$-plane.
}

Before continuing, it will help to reformulate the ordered induced matching language of Loh \cite[Lemma 2.2 and Observation 2.2]{Loh}.
Note that $x$-slices and $y$-slices of $S\belongs \ZZ^3$ correspond to the rows and columns of the $xy$-grid view, respectively, while the $z$-slices correspond to labeled squares with a given label.

\obs{
\label{OBS:basic-grid-observations}

Let $S\belongs \ZZ^3$ be a slice-increasing (or ordered) set.
Then in the $xy$-grid view of $S$, any fixed row or column has increasing labels, and any fixed label $z$ appears in coordinate-wise increasing squares $(x_1,y_1)<\dots<(x_m,y_m)$.
Furthermore, the ``$z$-corners'' $(x_i,y_j)$ and $(x_j,y_i)$ are empty for each pair $i<j$ \cite[Lemma 2.2]{Loh}.
}

\subsection{Ramsey for ordered surfaces, and connection from \Szabo--Tardos}
\label{SUBSEC:Ramsey-and-Szabo-Tardos}

We first show that any ``flat'' subset of a \emph{constrained} slice-increasing subset is small.
With care, one could reformulate this section in terms of \emph{ordered surfaces}: see Section \ref{SUBSEC:ordered-surfaces-appendix}.

\prop{
\label{PROP:primary-ordered-surfaces-are-small}

If a slice-increasing set $S\belongs [n]^3$ has a subset $T$ \emph{avoiding at least one of the difference-types $(+,+,\leqslant\!0), (+,\leqslant\!0,+), (\leqslant\!0,+,+)$}, then $\card{T} \le n$.
}

\pf{
Suppose $T$ avoids $(+,+,\leqslant\!0)$.
Since $T$ is slice-increasing, it contains at most a single point on every $z$-slice, so $\card{T} \le n$.
}

As observed in Section \ref{SUBSUBSEC:K-flatness}, a naive application of Theorem \ref{THM:main-theorem} only shows that a \K-free geometric tournament $\mcal{T} = \Color(T)$ with $T\belongs [n]^3$ has $\card{T} \le n^{3/2}$.
But we can do much better, not only for ordered sets, but for any slice-increasing set.

\prop{
\label{PROP:linear-bound-on-K-free-slice-increasing-sets}

If a slice-increasing set $S\belongs [n]^3$ has a subset $T$ \emph{avoiding difference-type $(+,+,+)$}, then $\card{T} \le 3n$.
}

\pf{
\WLOG{} $S = T$, and use the $xy$-grid view.
Purge rows with at most 1 point to get a new slice-increasing set $S'$ of size at least $\card{S} - n$.
If column $x$ of $S'$ has 3 points $P_1 = (x,y_1,z_1)$, $P_2 = (x,y_2,z_2)$, $P_3 = (x,y_3,z_3)$ in increasing order, then row $y_2$ cannot contain any other points of $S$: a point $(x',y_2,z')$ with $z' > z_2$ would create a $(+,+,+)$-type difference directed from $P_1$, while a point $(x',y_2,z')$ with $z' < z_2$ would create a $(+,+,+)$-type difference directed towards $P_3$.
Thus every column of $S'$ in fact has at most 2 points, so $\card{S'} \le 2n$ and we conclude $\card{S} \le \card{S'}+n \le 3n$.
}

\subsubsection{\Szabo--Tardos}

In \cite{SzaboTardos}, \Szabo and Tardos asked for the asymptotics of $m(N,2)$, the largest number $M$ such that any set $S\belongs \RR^3$ of size $N$ has a subset $T$ of size $M$ that avoids at least one of the four \emph{strict} difference-types $(+,+,+), (+,+,-), (+,-,+), (-,+,+)$.

\rmk{
Strictly speaking, they also have a requirement that all of the $x$-coordinates are distinct.
However, by slightly perturbing the points of $S$ into coordinate-general position, one can only make it harder to find such $T$, because \emph{strict} difference-types are stable under small perturbation.
So both versions of the problem are equivalent to the version with $S$ required to be in coordinate-general position.
}

\ex{[\Szabo--Tardos {\cite[Theorem 5, $(d,i) = (2,2)$ case]{SzaboTardos}}]

There are arbitrarily large $N$ with $m(N,2) \le CN^{5/8}$, for some absolute constant $C>0$.
The examples $S$ achieving this bound can be taken in coordinate-general position.
}

In particular, it is \emph{not} true that $m(N,2) \ge N^{2/3}$ for all $N$, as one might initially hope.
However, as \Szabo and Tardos suggest \cite[Section 4, Remark 3]{SzaboTardos}, it still seems likely that $m(N,2)$ is substantially larger than the trivial bound $N^{1/2}$.
We now show that such a nontrivial lower bound would give nontrivial upper bounds on slice-increasing sets.

\thm{
\label{THM:Szabo-Tardos-transfers-to-slice-increasing-bounds}

Let $S\belongs [n]^3$ be a slice-increasing set.
Then $m(\card{S},2) \le 3n$.
In particular, if there exists $\alpha > 1/2$ such that $m(N,2) \ge N^{\alpha}$ for all $N$, then $\card{S} \le n^{1/\alpha}$.
}

\pf{
Perturb $S$ to $S' = \phi(S)$ via a real linear map $\phi\maps (x,y,z)\mapsto (x-\eps z, y - \eps x, z - \eps y)$ for some \emph{positive} $\eps < n^{-1}/10$ such that $\phi$ is invertible.
Observe that $\phi$ converts difference-types $(+,+,0)$ into type $(+,+,-)$.
Take a subset $T'$ of $S'$ of size $m(\card{S},2)$, with inverse $T$, such that $T'$ avoids one of the four possible strict difference-types.
\begin{itemize}
\item If $T'$ avoids $(+,+,+)$, then so does $T$.
So in this case, Proposition \ref{PROP:linear-bound-on-K-free-slice-increasing-sets} gives $m(\card{S},2) = \card{T'} = \card{T} \le 3n$, as desired.

\item Otherwise, by cyclic symmetry, \Wlog{} suppose $T'$ avoids $(+,+,-)$.
We claim that $T$ avoids $(+,+,\leqslant\!0)$.
Indeed, suppose $(x,y,z)$ and $(u,v,w)$ lie in $T$ such that $u>x$; $v>y$; and $w \le z$.
Then after perturbation, $u' > x'$ and $v' > y'$ still.
If $w < z$, then similarly $w' < z'$; but $w' < z'$ even if $w=z$, by the observation that $\phi$ converts $(+,+,0)$ differences into $(+,+,-)$ differences.
So regardless, the perturbed difference is type $(+,+,-)$, contradicting the assumption on $T'$.
So indeed $T$ avoids $(+,+,\leqslant\!0)$, and $m(\card{S},2) = \card{T'} = \card{T} \le n$ by Proposition \ref{PROP:primary-ordered-surfaces-are-small}.
\end{itemize}

This shows $m(\card{S},2) \le 3n$.
Now suppose $m(N,2) \ge N^\alpha$ uniformly for some exponent $\alpha>1/2$.
Then $\card{S}^\alpha \le 3n$, so $\card{S} \le 3^{1/\alpha}n^{1/\alpha}$ for all $n$ and $S$.
To remove the constant in front, we consider lexicographic ``powers'' of $S$, as in Remark \ref{RMK:tensor-power-trick}.
}

\rmk{
We give one possible heuristic for why the \Szabo--Tardos problem cannot get the expected $n^{3/2}$ bound on the slice-increasing sets problem.
Consider the perturbation $\phi$ used in the proof of Theorem \ref{THM:Szabo-Tardos-transfers-to-slice-increasing-bounds}.
If $S$ is \emph{ordered}, then $\phi(S)$ has the same ordering---in fact we have an equality of geometric tournaments, $\Color(\phi(S)) = \Color(S)$.
If we properly extend $\Color$ to slice-increasing sets, then we instead have an equality of \emph{non-transitive} tournaments.
If one could prove that every $N$-vertex non-transitive 3-colored tournament has a simple directed 1-color-avoiding path of length $N^{2/3}$, then the slice-increasing problem would be completely resolved.
However, if one applies \Szabo--Tardos to $\phi(S)$, then one is looking for 1-color-avoiding \emph{cliques} in the \emph{non-transitive} 3-colored tournament $\Color(\phi(S)) = \Color(S)$.
One can show every cyclic triangle in $\Color(S)$ is rainbow with edges $\R,\G,\B$ in some order, so 1-color-avoiding cliques in $\Color(S)$ are also 1-color-avoiding paths.
But a priori, 1-color-avoiding paths need not be 1-color-avoiding cliques.
(It may be that they are when $\card{S}$ is far from general position, but that would require nontrivial proof.)
}

\subsubsection{Edge-type and difference-octant counts}


The previous remark suggests that for the slice-increasing problem itself, it may help to do one of the following.
\begin{itemize}
\item From the Ramsey perspective, try not to find 1-color-avoiding cliques, but something closer to paths.
We do not know of a simple geometric way to think about paths without going through the transformation machinery of Section \ref{SEC:transformations}.
Furthermore, some additional complications arise for non-transitive tournaments, such as the need to avoid cycles when discussing paths.

\item Use the fact that when $\card{S} \gg n^{3/2}$, the set $S$ is \emph{far from general position}---the coordinate-slices in a given orientation contain at least $\card{S}/n\gg n^{1/2}$ points on average.
As a first example, one can extend the proof of Proposition \ref{PROP:linear-bound-on-K-free-slice-increasing-sets} to \emph{count} the number of $(+,+,+)$-type differences (say) appearing in a large slice-increasing set.
For more on the pair-counting approach, see Section \ref{SUBSEC:edge-counts-appendix}.
\end{itemize}

\subsection{Sum of squares of slice-counts}
\label{SUBSEC:L^2-slice-counts-approach}

By Observation \ref{OBS:basic-grid-observations}, if some $z$-slice has $m$ points, then in the $xy$-grid view, the corresponding $z$-labels occupy the main diagonal of an otherwise empty $m\times m$ grid.
But the known tight examples $S\belongs [n]^3$ for Question \ref{QUES:bounding-slice-increasing-sets} have the same number of points $n^{1/2}$ on each slice, so on average, these $m\times m$ grids of $m^2 = n$ points will cover each square of $[n]^2$ exactly once.
This motivates the following $L^2$ questions.
Positive progress on any of them would lead to positive progress on the original slice-increasing questions, by Cauchy--Schwarz.

\ques{
Let $S\belongs [n]^3$ be a slice-increasing set.
For $i\in [n]$, let $a_i$ denote the size of the $z$-slice $\set{z = i}\cap S$.
Is $\sum_i a_i^2$ always at most $n^{2.99}$?
}

\ex{
Whenever $\sum_i a_i^2 \le n^2$, one has $\sum_i a_i \le n^{3/2}$ by Cauchy--Schwarz.
However, there are arbitrarily large slice-increasing examples with $\sum_i a_i^2 = n^\alpha$, for some \emph{exponent} $\alpha > 2$.
To see this, it suffices---by the lexicographic ``powers'' trick in Remark \ref{RMK:tensor-power-trick}---to find a single slice-increasing example with $\sum_i a_i^2 \ge n^2 + 1$.
The smallest we have been able to find is the following set $S\belongs [6]\times [6]\times [4]$, shown in the $xy$-grid view.
\mathd{
\begin{tabular}{ | l | l | l | l | l | l | }
  \hline
    &   &   & 2 &   & 4\\ \hline
    &   &   &   &   & 1\\ \hline
    &   &   & 1 & 4 &  \\ \hline
  2 &   & 4 &   &   &  \\ \hline
    &   & 1 &   & 3 &  \\ \hline
  1 & 4 &   &   &   &  \\ \hline
\end{tabular}
}
Here the list $(a_1,a_2,a_3,a_4) = (4,2,1,4)$ has sum of squares $37 \ge 6^2 + 1$.
This example was constructed by ``doubling'' a ``tight'' example for $n=3$, and then inserting in one more number in the bottom-right fourth of the grid.
}

We have not been able to find similar ordered examples, due to the observation below.
So the following two questions are still open.
A positive answer to either would completely resolve the original ordered problems.

\ques{
Let $S\belongs [n]^3$ be an ordered set.
For $i\in [n]$, let $a_i$ denote the size of the $z$-slice $\set{z = i}\cap S$.
Is $\sum_i a_i^2$ always at most $n^2$?
}

\ques{
Let $S = \Record(\mcal{T})$ for a canonical tournament $\mcal{T}$.
For $i\in [\ell(\GBK)]$, let $a_i$ denote the size of the $z$-slice $\set{z = i}\cap S$.
Is $\sum_i a_i^2$ always at most $\ell(\RGK)\cdot\ell(\RBK)$?
}

\obs{

Consider the $xy$-grid view of an ordered set $S\belongs \ZZ^3$.
Suppose for some $a<b$, a pair of $a$-labeled squares share a ``corner'' with a pair of $b$-labeled squares.
Then by Observation \ref{OBS:basic-grid-observations} they form a $3\times3$ sub-grid as follows:
\mathd{
\begin{tabular}{ | l | l | l | }
  \hline
    &   & b \\ \hline
    &   & {\color{red}a} \\ \hline
  a & {\color{red}b} & \textup{x}  \\ \hline
\end{tabular}
\qquad\textup{or}\qquad
\begin{tabular}{ | l | l | l | }
  \hline
  \textup{x} & {\color{red}a} & b \\ \hline
  {\color{red}b} &   &   \\ \hline
  a &   &   \\ \hline
\end{tabular}
}
(here the \textup{x} marks the overlapping corner).
Furthermore, this $3\times3$ sub-grid is empty except for the aforementioned $a,b,a,b$ labels.
}

\pf{
\WLOG{} assume the \textup{x} is on the lower-right corner of the sub-grid.
Since $S$ is slice-increasing, Observation \ref{OBS:basic-grid-observations} already rules out all but the upper-left corner of the sub-grid.
A slice-increasing set might have a label $c$ in that square, forcing $a<c<b$.
But such a $c$-square would violate the ordering condition, Definition \ref{DEF:ordered-triples}, with the red squares: specifically, the $b$ on the bottom-most row, and the $a$ on the right-most column.
}

This intuitively seems powerful, but it is ``local'' so we have not figured out to use it well.

\subsubsection{Edge-type and difference-octant counts}

If $\sum_i a_i^2$ is large, then as a first ``global'' step, there must be many \R-edges and \K-edges in $\Color(S)$: see Section \ref{SUBSUBSEC:L^2-large-implies-edge-lower-bound}.

\subsubsection{Representation theory}

These questions bear some resemblance to the following foundational result from representation theory.
It would be interesting if it could be applied to the problem at hand, since to our knowledge most applications of representation theory to extremal combinatorics are instead based on character theory.

\thm{[{\cite[Corollary 3.5.5 in book; Corollary 2.13 in arXiv v5]{Etingof}}]
Let $A$ be a finite-dimensional unital associative algebra over an algebraically closed field $k$.
Then $A$ has finitely many irreducible representations $V_i$ up to isomorphism, and $\sum_i \dim_k(V_i)^2 \le \dim_k(A)$.
}

\subsection{Overlap with joints problem}
\label{SUBSEC:joints-problem-inspiration}

Example \ref{EX:standard-slice-increasing-tight-example-sumsets} gives slice-increasing sets $S\belongs[n]^3$ with $n^{3/2}$ points, defined by certain coordinate-wise strictly increasing functions $f,g,h$.
Suppose $f,g,h$ are affine maps, e.g. $(r,s)\mapsto n^{1/2}r+s$ in the most standard example.
Then each of the $3n$ coordinate-slices of $S$ consists of collinear points lying on some \emph{increasing} line $\ell$.
Let $L$ be the set of these $3n$ lines.
Then each point $P_0 = (x_0,y_0,z_0)$ of $S$ lies on exactly three lines $\ell_1,\ell_2,\ell_3$ of $L$, where $\ell_1$ (resp. $\ell_2$; $\ell_3$) denotes the line of $L$ contained in the slice $x = x_0$ (resp. $y = y_0$; $z = z_0$).
Since $\ell_1,\ell_2,\ell_3$ are \emph{increasing} on their respective slices, one easily checks that they are \emph{non-coplanar}.

\obs{
In this example, each of the $n^{3/2}$ points of $S$ is a \emph{joint} of $L$: an intersection of three \emph{non-coplanar} lines.
But $\card{L} = 3n \sim n$, so up to an absolute constant, $L$ is a tight example for the joints problem \cite{GuthKatz}.
}

\rmk{
One can say more, at the very least when $f=g=h\maps(r,s)\mapsto n^{1/2} r+s$.
Here there are three canonical ways to partition $S$ into $n^{1/2}$ planes with $n$ points each, and any two of these $3n^{1/2}$ planes with different ``orientations'' (cf. Section \ref{SUBSEC:ordered-surfaces-appendix}) intersect in one of the lines of $L$.
}


Of course, for slice-increasing sets $S\belongs [n]^3$ in general, the $3n$ slices will not be lines, but instead \emph{monotone discrete curves}.
So the polynomial method of \cite{GuthKatz} (or subsequent simplifications) may not directly extend.
However, one may be able to leverage the special positions of these monotone curves.
Also, perhaps by using lexicographic products (Definition \ref{DEF:lex-product-triples}), one could reduce to a more polynomial-like situation.

\section*{Acknowledgements}

The first two authors thank their mentor, Ben Yang, for suggesting the problem, for suggesting the incidence geometry connections in Sections \ref{SUBSEC:L^2-slice-counts-approach} and \ref{SUBSEC:joints-problem-inspiration}, and for countless invigorating discussions throughout the MIT SPUR 2016 program.
They also thank Profs. David Jerison and Ankur Moitra for their support and advice, the SPUR program for the opportunity to conduct this research project, and Dr. Slava Gerovitch for directing the program.
The second author also thanks Po-Shen Loh for helpful lectures and classes throughout the years at the Math Olympiad Summer Program.

\appendix
\section{Python code, and example of a non-Gallai canonical tournament}
\label{SEC:Python-code}


We starting by linking to the following three Python programs, in case they are helpful.\footnote{Actual Python files, together with further examples, also accompany the arXiv version of this paper.}
\begin{enumerate}
\item \url{https://repl.it/Cc7k/61} generates random tournaments, implements the $\Color$, $\Record$, and $\Dual$ maps, and has a function checking whether an arbitrary \RGBK-tournament is undirected-Gallai (Definition \ref{DEF:undirected-Gallai-RGBK}).

\item Given an \RGBK-tournament of $\mcal{T}$, \url{https://repl.it/Ccyp/25} tests whether random blowups of $\mcal{T}$ satisfy Question \ref{QUES:Ramsey-for-RGBK-tournaments}.

\item \url{https://repl.it/Cci2/103} randomly tests Questions \ref{QUES:bounding-unequal-slice-increasing-sets} and \ref{QUES:original-triples-problem}.
\end{enumerate}

\begin{figure}[ht]
\includegraphics[scale=0.4]{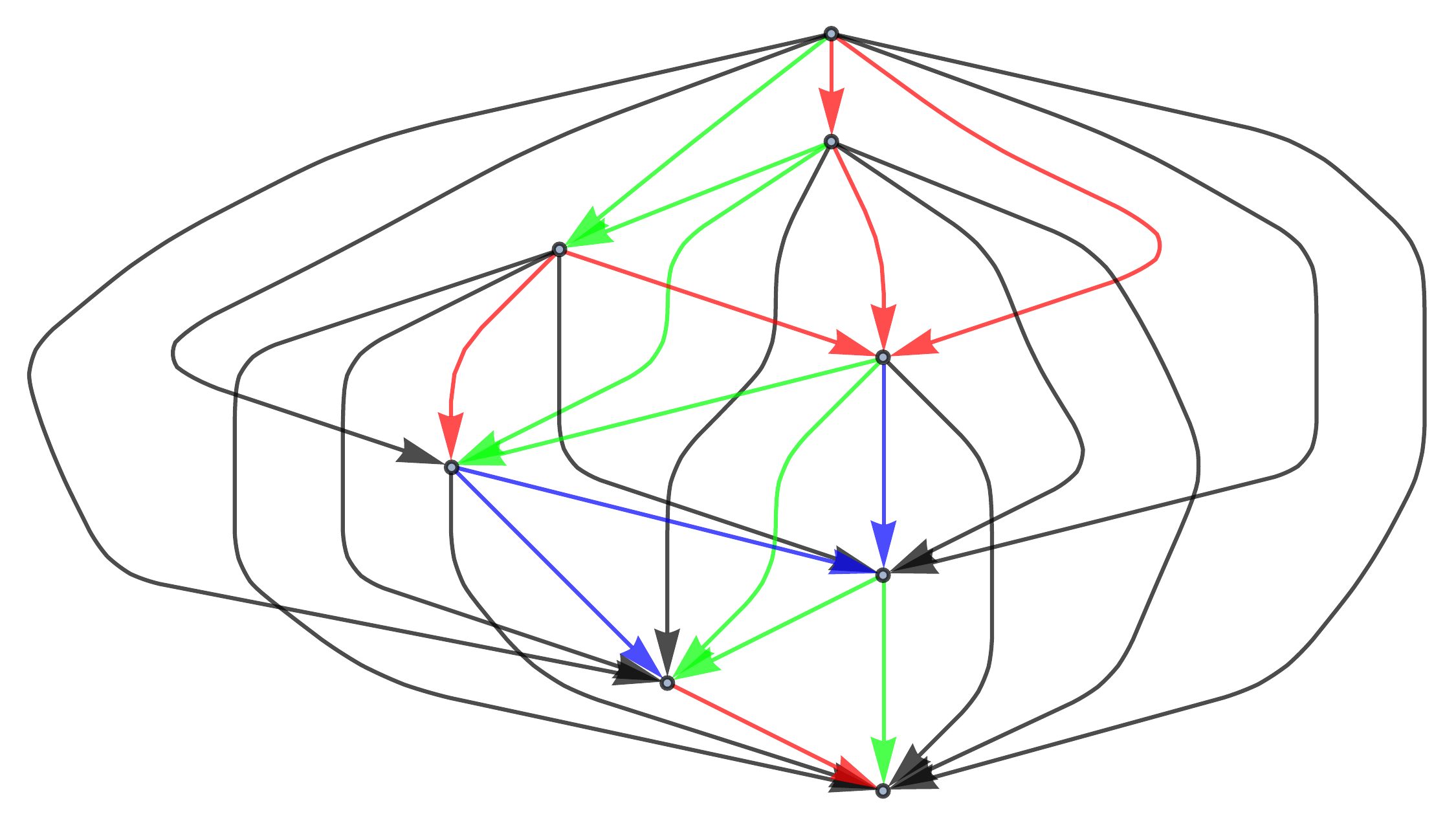}
\end{figure}

By random search we were able to find a canonical \RGBK-tournament $\mcal{T} = \Color(S)$ on $N=8$ vertices that is \emph{not} undirected-Gallai; the visualization here was created with Mathematica.
The ordered set $S = \Record(\mcal{T})$ is
\mathd{
\set{(1,1,1),(2,2,1),(3,1,2),(4,3,1),(5,2,2),(4,4,3),(5,3,4),(6,4,4)}.
}
To check that $\mcal{T}$ is canonical, one can verify the properties of $S$ required by the characterization in Theorem \ref{THM:classification-of-canonical-tournaments}.
To check that $\mcal{T}$ does not have an undirected, \K-blind Gallai decomposition, one notes that the blocks in a valid base decomposition must be unions of the connected components in whatever primary color the base graph avoids.



\section{Tournaments fixed under \texpdf{$\Color\circ\Record$}{Color-Record}}
\label{SEC:structure-of-Color-Record-fixed-points}

In this section we study \RGBK-tournaments in the image of $\Color\circ\Record$, which turn out to coincide with the fixed points of $\Color\circ\Record$.
The main result is the following classification, which we use in Appendix \ref{SEC:structure-of-canonical-tournaments} to classify canonical tournaments.

\subsection{Classification}

\thm{
\label{THM:classification-of-Color-Record-fixed-points}

Let $\mcal{T}$ be an \RGBK-tournament.
The following are equivalent.
\begin{enumerate}
\item $\mcal{T}$ is fixed under $\Color\circ\Record$.

\item $\mcal{T}$ lies in the image of $\Color\circ\Record$.

\item Tournament $\mcal{T}$ is geometric.
For each color class $\C\in \set{\RGK,\RBK,\GBK}$, the \C-stratification (Definition \ref{DEF:C-stratification}) of $\mcal{T}$ is geometric (Definition \ref{DEF:stratification-is-geometric}), and the \C-minimal vertices (Definition \ref{DEF:C-minimal-vertices}) form an increasing \C-colored path.

\item $\mcal{T} = \Color(S)$ for some \emph{reduced set} $S\belongs \ZZ_{>0}^3$: an ordered set of triples such that for each coordinate $c\in\set{x,y,z}$ and $\bd{p}\in S$ with $c$-coordinate $i\ge 2$, there exists a point $\bd{q}\in S$ with $c$-coordinate exactly $i-1$, such that $\bd{p} - \bd{q}$ has at least two positive coordinates.
\end{enumerate}
Furthermore, if $\mcal{T}$ satisfies any of these equivalent conditions, then the only $S$ that works in (4) is $S = \Record(\mcal{T})$.
}

\pf{[Proof of equivalences]
(1) clearly implies (2), while (2) implies (3) by Propositions \ref{PROP:stratification-class-1-properties} (stratification properties) and \ref{PROP:reduced-graph-structure} (reduced graph structure) below.

Now take a tournament $\mcal{T}$ satisfying (3).
Define $S \defeq \Record(\mcal{T})$.
We first show that $S$ is a reduced set.
\WLOG{} let $c = x$.
Consider the \RGK-stratification $V_1\sqcup\dots\sqcup V_m$ of $\mcal{T}$.
By definition, $x_i = j$ (i.e. the longest \RGK-path of $\mcal{T}$ ending at $v_i$ has length $j$) if and only if $v_i\in V_j$.
So if $(x_k,y_k,z_k) \in S$ with $x_k = i \ge 2$, then by considering $v^0_{i-1} < v^0_i \le v_k$ (where it is important that the stratification is geometric, and that the \RGK-minimal vertices form a \C-colored path), we see that the triple $(x',y',z')\in S$ corresponding to $v^0_{i-1}$ has $x$-coordinate $i-1$, and is strictly less than $(x_k,y_k,z_k)$ in at least two coordinates since $v^0_{i-1} < v_k$ in $\mcal{T}$.
So $S$ is reduced.
Also, if $v_i < v_j$ in $\mcal{T}$ is \B-colored (say), then by the geometric \RGK-stratification, $x_i$ is \emph{not} strictly less than $x_j$, so $v_i \to v_j$ cannot become \K-colored under $\Color\circ\Record$.
By Proposition \ref{PROP:Color-Record-RGK-stability}, we conclude that $\mcal{T} = \Color(\Record(\mcal{T})) = \Color(S)$.
So (3) implies (4) and (1).

Finally, take a tournament $\mcal{T} = \Color(S)$ satisfying (4).
To finish, it suffices to show that $\mcal{T}$ satisfies (2).
In fact, we claim that $S = \Record(\mcal{T})$.
Take the $k$th triple $(x_k,y_k,z_k) \in S$.
Then the longest \RGK-path in $\mcal{T}$ ending at the $k$th vertex has length not only \emph{at most} $x_k$ (since the $x$-coordinates of $S$ are all positive integers, and an \RGK-edge requires an $x$-increase in the triples of $S$), but \emph{at least} $x_k$, by the $x$-traceback property in the definition of reduced sets.
Thus $S = \Record(\mcal{T})$, as desired.
}

\subsection{Transitivity, and paths vs. cliques}

For the reader's convenience, we start with an immediate corollary of Proposition \ref{PROP:image-of-Color-is-blah-transitive}.

\cor{
\label{COR:blah-transitivity-stated-for-class-1-RGBK}

An \RGBK-tournament in the image of $\Color\circ\Record$ is transitive in each of the color classes \R, \G, \B, \K, \RK, \GK, \BK, \RGK, \RBK, \GBK.
As a corollary, the notions of directed path and (undirected) clique coincide for each of these color classes.
}

\subsection{Stability}

\prop{[$\RGK$-stability]
\label{PROP:Color-Record-RGK-stability}

Applying transformation $\Color\circ \Record$ to an \RGBK-tournament has the effect of possibly changing some edges to color \K, but leaving all other edges unchanged.
The same holds for the dual conjugate $\Dual\circ\Color\circ \Record\circ\Dual$.
}

\pf{[Proof for original transformation]
Let $\mathcal T$ be our original transitive tournament and $\mathcal T'$ be the result of the transformation.
Let $v<v'$ be two vertices of $\mathcal T$ such that $vv'$ is \XK-colored for some $\X\in \set{\R,\G,\B}$, say $\X = \B$.
Let $v$ correspond to triple $(x,y,z)$ and $v'$ to $(x',y',z')$ under transformation $\Record$.
We claim that $y<y'$ and $z<z'$.
This is simply because any \RBK-colored path of length $y$ ending at $v$ can be extended by $vv'$ to an \RBK-colored path of length $y+1$ ending at $v'$ and similarly for \GBK-colored paths.
Therefore the edge $(x,y,z)\to(x',y',z')$ has orientation $(\leqslant\!0,+,+)$ or $(+,+,+)$ so in $\mathcal T'$, the edge $vv'$ is \BK-colored.
Thus if an edge is color \XK in $\mcal{T}$ for some $\X\in \set{\R,\G,\B}$, then it is still color \XK in $\mcal{T}'$.
So an edge of color \RGB can only change to \K (or stay the same), and an edge of color \K cannot change.
}

\pf{[Proof for dual transformation]

Note that $\Dual$ merely reverses edge directions without changing colors.
So the statement for $\Color\circ\Record$ respects $\Dual$.
}

\subsection{Stratification}

\defn{
\label{DEF:C-stratification}
Given an \RGBK-tournament $\mcal{T}$, define the \emph{\C-stratification} as the vertex-partition $V_1\sqcup \dots \sqcup V_m$ where $v\in V_\ell$ if and only if the longest \C-colored path ending at $v$ has length $\ell$.
For example, if $\C = \RGK$, then $m = \ell(\RGK)$.
}

\defn{
\label{DEF:stratification-is-geometric}

If $\C\in\set{\RGK,\RBK,\GBK}$, call the \C-stratification $V_1\sqcup \dots \sqcup V_m$ \emph{geometric} if for $u\in V_i$ and $v\in V_j$ with $u<v$, the edge $u\to v$ is \C-colored if and only if $i<j$.
}

\prop{[Traceback]
\label{PROP:stratification-traceback}

Consider the \C-stratification $V_1\sqcup \dots \sqcup V_m$ of an arbitrary \RGBK-tournament $\mcal{T}$.
Then for each $v\in V_i$ with $i\ge2$, there is a vertex $u\in V_{i-1}$ with $u<v$ and \C-colored edge $u\to v$.
}

\pf{
Suppose $v\in V_i$ for $i>1$.
Then there exists an \C-colored path of length $i$ that ends at $v$.
Let $u$ be the second-to-last vertex in this path.
Clearly $u\in V_{i-1}$.
}

\prop{
\label{PROP:stratification-class-0-properties}

Consider the \C-stratification $V_1\sqcup \dots \sqcup V_m$ of an arbitrary \RGBK-tournament $\mcal{T}$.
Then $i<j$ for every \C-edge $u\to v$ with $u\in V_i$ and $v\in V_j$.
}

\pf{
If $u\to v$ is \C-colored, then any \C-colored path ending at $u$ can be extended by $uv$, so $j > i$.
}

Sometimes we can say more.

\prop{
\label{PROP:stratification-class-1-properties}

A tournament in the image of $\Color\circ\Record$ has geometric \C-stratification for each $\C\in\set{\RGK,\RBK,\GBK}$.
}

\pf{
\WLOG{} $\C = \RGK$.
Consider the \RGK-stratification $V'_1\sqcup \dots \sqcup V'_m$ of a tournament $\mcal{T}' = (\Color\circ\Record)(\mcal{T})$.
Proposition \ref{PROP:stratification-class-0-properties} (applied to $\mcal{T}'$, not $\mcal{T}$) says that color \RGK implies $i<j$, or equivalently that $i\ge j$ implies color \B.
It remains to prove the converse, i.e. that when $i < j$, every edge $u\to v$ from $u\in V'_i$ to $v\in V'_j$ is \RGK-colored.

For clarity, say $\mcal{T}'$ has vertex-list $v_1,\ldots, v_N$, while $\mathcal T$ has vertex-list $w_1,\dots,w_N$.
Let $(x_i,y_i,z_i)$ be the triple that corresponds to vertex $w_i$ under transformation $\Record$.
We claim that $v_i$ lies in $V'_\ell$ if and only if $\ell = x_i$, or in other words that the length $\ell$ of the longest \RGK-path in $\mcal{T}'$ ending at $v_i$ equals the length $x_i$ of the longest \RGK-path in $\mcal{T}$ ending at $w_i$.

\begin{itemize}
\item Let $v_{i_1},v_{i_2},\ldots, v_{i_j}=v_i$ be an \RGK-path in $\mathcal T'$ ending at $v_i$.
Then $x_{i_1}<x_{i_2}<\dots<x_{i_j}=x_i$.
These $x$-coordinates are positive integers, so $x_i\geq j$, and $x_i \ge \ell$.
(Cf. proof of Proposition \ref{PROP:reduction-by-Record}.)

\item Conversely, consider a length-$j$ \RGK-path in $\mathcal T$ ending at $w_i$.
By Proposition \ref{PROP:Color-Record-RGK-stability}, \RGK is preserved under $\Color\circ \Record$, so such an \RGK-path becomes an \RGK-path in $\mathcal T'$ ending at $v_i$.
Thus $\ell \ge j$, meaning $\ell \ge x_i$.
\end{itemize}
This proves the claim.
Finally, suppose $v_i\in V'_\ell$ and $v_{i'}\in V'_{\ell'}$ with $v_i < v_{i'}$ and $\ell < \ell'$.
Then $x_i = \ell < \ell' = x_{i'}$ by the claim, so by definition of $\mcal{T}'$, edge $v_i\to v_{i'}$ is \RGK-colored, as desired.
}

\rmk{
Intuitively, we obtain geometric properties from $\Color$, combinatorial properties from $\Record$, and a compatible mixture from their composition.
The proof also shows, up to definition-chasing, that $\Record\circ\Color\circ\Record = \Record$, so $\Color\circ\Record$ is idempotent.
}

\subsection{\C-minimal vertices and \C-reduced graph}

\defn{
\label{DEF:C-minimal-vertices}

For a tournament $\mcal{T}$ with \C-stratification $V_1\sqcup \dots \sqcup V_m$, define $v^0_i\in V_i$ to be the \emph{smallest} vertex in set $V_i$ and let the \emph{\C-reduced graph} $\mcal{T}^0$ be the transitive sub-tournament induced by these \emph{\C-minimal vertices} $v^0_1,\ldots,v^0_m$.
}

\prop{
\label{PROP:reduced-graph-structure}

Suppose $\mcal{T}$ lies in the image of $\Color\circ\Record$.
Then the \RGK-reduced graph $\mcal{T}^0$ is \RGK-colored with $v^0_1<v^0_2<\cdots<v^0_m$, and transitive in each of \RK, \GK, \R, \G, \K.
Analogous statements hold for the \RBK-reduced graph and the \GBK-reduced graph.
}

\pf{
Proposition \ref{PROP:stratification-traceback} (traceback) implies $v^0_m > v^0_{m-1} > \dots > v^0_1$.
In particular, $v^0_i<v^0_j$, so edge $v^0_i v^0_j$ must be \RGK-colored by Proposition \ref{PROP:stratification-class-1-properties}.
For transitivity see Corollary \ref{COR:blah-transitivity-stated-for-class-1-RGBK}.
}

In the proof of Theorem \ref{THM:classification-of-canonical-tournaments} we will use the following stratification-free characterization of minimal vertices.

\prop{
\label{PROP:stratification-free-characterize-minimal-vertex}

Suppose $\mcal{T}$ lies in the image of $\Color\circ\Record$, and fix a color class $\C\in\set{\RGK,\RBK,\GBK}$.
Then $v^0_1 = 1$ is the smallest vertex in $\mcal{T}$, while for $i\ge2$, the $i$th \C-minimal vertex $v^0_i$ is the smallest vertex $v > v^0_{i-1}$ with a \C-colored edge from $v^0_{i-1}$.
}

\pf{
For $i=1$, note that the first vertex clearly lies in $V_1$, so $v^0_1 = 1$.
Now suppose $i\ge2$.
Then by Proposition \ref{PROP:stratification-class-1-properties}, a vertex $v\in \mcal{T}$ has \C-colored edge \emph{from} $v^0_{i-1}$ if and only if $v\in V_i\cup\dots\cup V_m$; it is important here that $v^0_{i-1}$ is the smallest vertex in $V_{i-1}\cup V_i \cup \dots\cup V_m$, by Proposition \ref{PROP:reduced-graph-structure}.
But again by Proposition \ref{PROP:reduced-graph-structure}, the smallest vertex in $V_i\cup\dots\cup V_m$ is $v^0_i$, as desired.
}

\subsection{\B-connected components}

The following extends the interval-connectivity observation from Proposition-Definition \ref{PROPDEF:RGBK-interval-connected}.

\prop{
\label{PROP:basic-description-B-connected-components}

Suppose $\mcal{T}$ lies in the image of $\Color\circ\Record$.
Then the connected \B-components of $\mcal{T}$ are interval-connected.
Each \B-component is a union of \emph{consecutive} \RGK-strata $V_i$.
Furthermore, for $i<j$, \emph{if} the \RGK-strata $V_i$ and $V_j$ lie in \emph{different} \B-components, \emph{then} $V_i < V_j$ absolutely, i.e. $u < v$ for all $u\in V_i$ and $v\in V_j$.
}

\pf{
Proposition-Definition \ref{PROPDEF:RGBK-interval-connected} gives interval-connectivity, since $\mcal{T}$ is geometric.
Each \RGK-stratum is a \B-clique by the easy half of Proposition \ref{PROP:stratification-class-1-properties}, hence trivially \B-connected.
Now take indices $i<j<k$, and suppose that $V_i,V_k$ are in the same \B-component.
Since $v^0_i < v^0_j < v^0_k$ by Proposition \ref{PROP:reduced-graph-structure}, interval-connectivity implies that $v_0^j$ also lies in the same \B-component, so $V_j$ lies in the same \B-component.

Furthermore, for $i<j$, the \RGK-strata $V_i$ and $V_j$ have no \B-edges between them if and only if there are no pairs of vertices $u\in V_j$ and $v\in V_i$ with $u < v$.
(See Proposition \ref{PROP:stratification-class-1-properties}.)
This is in turn equivalent to having $V_i<V_j$ absolutely.
}

\cor{

The \B-components of a tournament $\mcal{T}$ in the image of $\Color\circ\Record$ are absolutely ordered, with all edges directed from smaller components to larger components.
}

\section{Canonical tournaments}
\label{SEC:structure-of-canonical-tournaments}


\subsection{Classification}

\thm{
\label{THM:classification-of-canonical-tournaments}

Let $\mcal{T}$ be an \RGBK-tournament.
The following are equivalent.
\begin{enumerate}
\item $\mcal{T}$ is canonical, i.e. fixed under $\Color\circ\Record$ and $\Dual\circ\Color\circ\Record\circ\Dual$.

\item Tournament $\mcal{T}$ is geometric.
For each color class $\C\in \set{\RGK,\RBK,\GBK}$, the \C-stratification (Definition \ref{DEF:C-stratification}) of $\mcal{T}$ is geometric (Definition \ref{DEF:stratification-is-geometric}), and the \C-minimal vertices (Definition \ref{DEF:C-minimal-vertices}) form an increasing \C-colored path, as do the \C-maximal vertices (Definition \ref{DEF:dual-stratification-and-maximal-vertices}).
The dual \C-strata (Definition \ref{DEF:dual-stratification-and-maximal-vertices}) are the usual \C-strata \emph{in reverse}.

\item $\mcal{T} = \Color(S)$ for some ordered set $S$ such that $S$ and $(1+\ell(\RGK),1+\ell(\RBK),1+\ell(\GBK)) - S$ are both reduced sets of triples belonging to $\ZZ_{>0}^3$ (see Theorem \ref{THM:classification-of-Color-Record-fixed-points}).
\end{enumerate}
Furthermore, if $\mcal{T}$ satisfies any of these equivalent conditions, then the only $S$ that works in (3) is $S = \Record(\mcal{T})$.
}

\pf{[Proof of equivalences]
First suppose (1) holds.
Both $\mcal{T}$ and $\Dual(\mcal{T})$ are fixed points of $\Color\circ\Record$, so by Theorem \ref{THM:classification-of-Color-Record-fixed-points} and Proposition \ref{PROP:canonical-max-vertices-properties}, (2) holds as well.
On the other hand, suppose (2) holds.
Then by Theorem \ref{THM:classification-of-Color-Record-fixed-points}, $\mcal{T}$ is fixed under $\Color\circ\Record$.
Since the dual \C-strata coincide with the usual \C-strata in reverse, the \C-stratification of $\Dual(\mcal{T})$ is geometric.
And since the \C-maximal vertices are increasing, the \C-minimal vertices of $\Dual(\mcal{T})$ are increasing as well.
So Theorem \ref{THM:classification-of-Color-Record-fixed-points} also shows that $\Dual(\mcal{T})$ is fixed under $\Color\circ\Record$, so (1) holds.

To finish, we prove the equivalence of (2) and (3).
If (2) and (1) hold, then Theorem \ref{THM:classification-of-Color-Record-fixed-points} first shows that $\mcal{T} = \Color(S)$ for $S \defeq \Record(\mcal{T})$, and $S$ is reduced.
One then checks that $(1+\ell(\RGK),1+\ell(\RBK),1+\ell(\GBK)) - S$ coincides with $\Record(\Dual(\mcal{T}))$, which is reduced by Theorem \ref{THM:classification-of-Color-Record-fixed-points}.
So (3) holds.
Conversely, if (3) holds, then one can check that $\Dual(\mcal{T}) = \Color((1+\ell(\RGK),1+\ell(\RBK),1+\ell(\GBK)) - S)$, so (2) holds by Theorem \ref{THM:classification-of-Color-Record-fixed-points} applied in reverse.
}

\defn{
\label{DEF:dual-stratification-and-maximal-vertices}

Given a color class \C and an \RGBK-tournament $\mcal{T}$, define the \emph{dual \C-strata} and \emph{\C-maximal vertices} as the strata and vertices of $\mcal{T}$ corresponding under $\Dual$ to the \C-strata and \C-minimal vertices of $\Dual(\mcal{T})$, respectively.
}

The term ``maximal'' is unambiguous by the following proposition.

\prop{
\label{PROP:canonical-max-vertices-properties}

If $\mcal{T}$ is canonical, then its \RGK-maximal vertices coincide with the maximal vertices in the usual \RGK-strata $V_1,\dots,V_m$, and $V_m,\dots,V_1$ are the dual \RGK-strata.
}

\pf{
We prove by strong induction on $i\ge1$ that the $i$th dual stratum is precisely $V_{m+1-i}$, and that the $i$th maximal \RGK-vertex is the largest vertex in $V_{m+1-i}$.
Suppose $i\ge1$ and assume the result holds up to $i-1$.
Take the maximal vertex $w_i$ in the $i$th dual stratum.
By the dual of Proposition \ref{PROP:reduced-graph-structure}, there exists a length $\ell(\RGK)+1-i = m+1-i$ \RGK-path ending at vertex $w_i$, so $w_i$ lies in $V_{m+1-i}\cup \dots\cup V_m$, hence in $V_{m+1-i}$ (as $V_m,\dots,V_{m+2-i}$ are the dual strata up to $i-1$, by the inductive hypothesis).
On the other hand, by the dual of Proposition \ref{PROP:stratification-free-characterize-minimal-vertex}, $w_i$ is
\begin{itemize}
\item the largest vertex $w$ in $\mcal{T} = V_1\cup\dots\cup V_m$ if $i=1$; and

\item the largest vertex $w$ with an \RGK-edge directed towards $w_{i-1}$ (the largest vertex of $V_{m+2-i}$, by the inductive hypothesis) if $i\ge2$.
\end{itemize}
Either way, by Proposition \ref{PROP:stratification-class-1-properties} we conclude that $w_i$ is the largest vertex $w\in V_1\cup\dots\cup V_{m+1-i}$, so every vertex $v\in V_1\cup\dots\cup V_{m-i}$ has an \RGK-edge directed towards $w_i$ by Proposition \ref{PROP:stratification-class-1-properties}.
So the $i$th dual stratum is a subset of $V_{m+1-i}$, while $V_{m+1-i}$ is a priori a subset of the union of the first $i$ dual strata.
But the first $i-1$ dual strata are $V_m,\dots,V_{m+2-i}$ by the inductive hypothesis.
So the $i$th dual stratum is precisely $V_{m+1-i}$, with largest vertex $w_i$, completing the induction.
}

\subsection{Why canonical?}

We now give two results suggesting that canonical tournaments are indeed ``canonical'' in some meaningful way.

\cor{
\label{COR:every-vertex-belongs-to-max-paths}

Fix $\C\in \set{\RGK,\RBK,\GBK}$.
Each vertex $v$ of a canonical tournament $\mcal{T}$ belongs to an \C-colored path of maximum vertex-length $\ell(\C)$.
In particular, under $\Record$, the last vertex of $\mcal{T}$ has coordinates $(\ell(\RGK),\ell(\RBK),\ell(\GBK))$.
}

\pf{
Use minimal vertices before $v$, and maximal vertices after $v$.
}

The following result contrasts canonical tournaments with general tournaments, whose edges can often be changed to \K without increasing longest 1-color-avoiding path lengths, by Proposition \ref{PROP:Color-Record-RGK-stability}.

\cor{[\K-saturation]
\label{COR:K-saturation-of-canonical}

If some edge $i\to j$ colored $\X\in\set{\R,\G,\B}$ in a canonical tournament $\mcal{T}$ is changed to \K, then the length of the longest \X-avoiding path strictly increases (and as usual, the other two weakly increase; in fact, they stay the same).
}

\pf{
Say $\X = \B$.
Then $j$ is either in the same, or an earlier, \RGK-stratum as $i$.
There is a maximum \RGK-path ending at $i$ using minimal vertices, and a maximum \RGK-path starting at $j$ using maximal vertices.
If the edge $i\to j$ is changed to \K, then these two paths can be concatenated to increase the length of the longest \RGK-path (in fact, by exactly 1 plus the difference of the stratum numbers of $j$ and $i$).
}

\subsection{Viewing canonical tournaments in grids}

\prop{

Let $\mcal{T}$ be canonical, and consider the $xy$-grid view (Definition \ref{DEF:projection-and-grid-view-of-triples}) of $S = \Record(\mcal{T})$.
Then for any $(x,y)\in [\ell(\RGK)]\times[\ell(\RBK)]$, there exists both
\begin{itemize}
\item a labeled square immediately above (i.e. in same column) or immediately to the right (i.e. in same row) of $(x,y)$, \emph{unless and only unless} $(x,y)$ is labeled and simultaneously maximal in its row and column; and

\item a labeled square immediately below (i.e. in same column) or immediately to the left (i.e. in same row) of $(x,y)$, \emph{unless and only unless} $(x,y)$ is labeled and simultaneously minimal in its row and column.
\end{itemize}
}

\pf{[Proof of first bullet point]
The only difficulty is addressing the case where $(x,y)$ is unlabeled.
In this case, suppose, for contradiction, that both of the following hold.
\begin{itemize}
\item There is no labeled square immediately to the right of $(x,y)$, or equivalently, the maximum in \emph{row} $y$ is $(u,y,a)$ for some $u < x$.

\item There is no labeled square immediately above $(x,y)$, or equivalently, the maximum in \emph{column} $x$ is $(x,v,b)$ for some $v < y$.
\end{itemize}
Now let $(u,y',a')$ be the maximum in \emph{column} $u$, so that $y'\ge y$ and $a'\ge a$.
Similarly, let $(x',v,b')$ be the maximum in \emph{row} $v$, so that $x' \ge x$ and $b' \ge b$.

The point is that $u<x$ and $y'\ge y > v$ together imply $a' < b$ by the required ordering of the \RGK-maximal vertices $(u,y',a')\to (x,v,b)$, while $v<y$ and $x' \ge x > u$ together imply $b' < a$ by the required ordering of the \RBK-maximal vertices $(x',v,b') \to (u,y,a)$.
But the former implies $a < b$ (as $a\le a'$) while the latter implies $b < a$ (as $b\le b'$), giving our desired contradiction.
}

\pf{[Proof of second bullet point]

This is entirely dual and only requires $\mcal{T}$ to be fixed under $\Color\circ\Record$.
}

\cor{
\label{COR:min-above,max-below}

Let $\mcal{T}$ be canonical, and consider the $xy$-grid view (Definition \ref{DEF:projection-and-grid-view-of-triples}) of $S = \Record(\mcal{T})$.
Then for any $(x,y)\in [\ell(\RGK)]\times[\ell(\RBK)]$ not in the $xy$-projection $S_z$, there exists a \emph{smallest} $z$-coordinate among squares immediately above or to the right of $(x,y)$, as well as a \emph{largest} $z$-coordinate among squares immediately below or to the left of $(x,y)$.
}

\rmk{
This is saying something about how row $y$ interacts with column $x$.
Can we do some counting along these lines?
Is there a useful three-dimensional fact along these lines?
}

\section{Slice-increasing sets}
\label{SEC:slice-increasing-sets}

\subsection{Ordered surfaces}
\label{SUBSEC:ordered-surfaces-appendix}

We start with the following tentative definition.

\defn{
A \emph{strict ordered surface} is the zero-locus $f(x,y,z) = 0$ of a coordinate-wise strictly monotone function $f\maps \ZZ^3 \to \RR$.
The \emph{orientation} type $\vec{c}\in\set{+,-}^3$ records the coordinate-wise directions of $f$, with a $+$ (resp. $-$) corresponding to a coordinate-wise increase (resp. decrease).
}

\defn{[Cf. \cite{SzaboTardos}]
Fix an orientation $\vec{c}\in \set{+,-}^3$.
Say a set of points \emph{completely avoids} the orientation $\vec{c}$ if there do not exist \emph{distinct} points $\bd{p},\bd{q}$ of the set with coordinate-wise difference $\bd{q} - \bd{p}$ \emph{weakly of sign type} $\vec{c}$, meaning of sign type $\set{\geqslant\!0,\leqslant\!0}^3$ where each \emph{strict} $+$ (resp. $-$) in $\vec{c}$ is replaced with a \emph{weak} $\geqslant\!0$ (resp. $\leqslant\!0$).
}

\thm{
\label{THM:characterization-of-strict-ordered-surface-points}

A finite set $T\belongs \ZZ^3$ is contained in a strict ordered surface with orientation $\vec{c}\in \set{+,-}^3$ if and only if $T$ completely avoids the orientation $\vec{c}$.
}

\pf{
By symmetry, assume $\vec{c} = (+,+,+)$.
First suppose $T \belongs\set{f = 0}$ for some coordinate-wise strictly increasing function $f\maps \ZZ^3 \to \RR$.
Suppose $\bd{p} = (x,y,z)$ and $\bd{q} = (x',y',z')$ are points of $T$ with $\bd{q}-\bd{p}$ weakly of sign type $\vec{c}$, i.e. $x'\ge x$; $y' \ge y$; and $z' \ge z$.
Then $0 = f(x',y',z') \ge f(x,y,z) = 0$.
In particular, equality holds, so strictness of $f$ forces $\bd{q} = \bd{p}$.
Thus $T$ completely avoids the orientation $\vec{c}$.

Conversely, suppose $T$ completely avoids the orientation $\vec{c} = (+,+,+)$.
We construct a $(+,+,+)$-surface through $T$ by induction on $\card{T}$, where the base case $\card{T}=1$ is obvious.
Now suppose $\card{T}\ge2$ and take a strict ordered surface $g(x,y,z) = 0$, with $g\maps \ZZ^3 \to \RR$ increasing in $x,y,z$, containing all of $T$ except possibly a point $P_0 = (x_0,y_0,z_0)$.
If $g(P_0)$ is zero, we are done.
Now suppose it is nonzero.
\begin{itemize}
\item If $g(P_0) < 0$, then we modify as follows to get the desired $f$: add $\abs{g(P_0)} = -g(P_0)$ to $g$ whenever $(x,y,z) \ge (x_0,y_0,z_0)$ coordinate-wise.

\item If $g(P_0) > 0$, then we instead modify as follows: subtract $\abs{g(P_0)} = g(P_0)$ from $g$ whenever $(x,y,z) \le (x_0,y_0,z_0)$ coordinate-wise.
\end{itemize}
In each case, $f|_{T\setminus P_0} = g|_{T\setminus P_0} = 0$ since $T$ completely avoids orientation $\vec{c}$.
So the resulting $f$ vanishes on $P_0$, still vanishes on $T\setminus P_0$, and still remains strictly increasing in each coordinate.
This completes the induction.
}

\cor{
\label{COR:surfaces-formulation-of-total-order}

A slice-increasing set $S\belongs \ZZ^3$ is ordered if and only if every three points of $S$ lie on a strict ordered surface with orientation $(+,+,-)$, $(+,-,+)$, or $(-,+,+)$.
}

\pf{
First, suppose $S$ is ordered and take three points $P_1 < P_2 < P_3$ in order.
By pigeonhole, \Wlog{} assume their $x$-coordinates are strictly increasing.
Then their $x$-coordinates are distinct, and whenever $x_j - x_i < 0$, either $y_j - y_i < 0$ or $z_j - z_i < 0$.
So these three points completely avoid the orientation $\vec{c} = (-,+,+)$, which by Theorem \ref{THM:characterization-of-strict-ordered-surface-points} gives a $\vec{c}$-type strict ordered surface through $P_1,P_2,P_3$.

Conversely, suppose $S$ is \emph{not} ordered.
Then the well-defined majority-comparable tournament on $S$ is \emph{not} acyclic, so there are three points $P_1,P_2,P_3$ forming a cycle.
By cyclicity, no edge can be weakly of sign type $(+,+,+)$.
Since $S$ is slice-increasing, each edge must instead be \emph{strictly} of sign type $(+,+,-)$, $(+,-,+)$, or $(-,+,+)$.
Again by cyclicity, no two edges can have the same sign type.
So $P_1,P_2,P_3$ do \emph{not} completely avoid any of the orientations among $(+,+,-)$, $(+,-,+)$, and $(-,+,+)$.
Thus by Theorem \ref{THM:characterization-of-strict-ordered-surface-points}, these three points do not lie on a strict ordered surface of any permissible orientation.
}

\ex{
The three points $(0,0,0)$, $(1,1,-2)$, and $(2,-3,1)$ are ordered.
Even though they lie on the plane $x+y+z = 0$---a strict ordered surface with orientation $(+,+,+)$ not listed---they also lie on a strict ordered surface of permissible orientation.
}

\subsection{Edge-counts}
\label{SUBSEC:edge-counts-appendix}

First, we note that the ``overlapping $z$-corners'' idea of Section \ref{SUBSEC:L^2-slice-counts-approach} and the ``L-shape''-counting method of Proposition \ref{PROP:linear-bound-on-K-free-slice-increasing-sets} extend separately to give the following edge-type lower bounds.
We refer to differences of type $(+,+,\leqslant\!0)$ as \R-edges, and so on, by slight extension of the $\Color$ definition for ordered sets (Definition \ref{DEF:triples-to-tournament}).

\thm{

Unless otherwise stated, let $S\belongs [n]^3$ be a slice-increasing set of size $\alpha n$.
\begin{enumerate}
\item An arbitrary set $S\belongs [n]^3$ has at least $\frac12\alpha(\alpha-1) n$ \R-edges (resp. \G-; \B-) coming from pairs of points on the same $z$-slice (resp. $y$-; $x$-).

\item If $\alpha \ge 3n^{1/2}$, then $S$ has at least $\frac1{16} \alpha^4$ \R-edges (resp. \G-; \B-) coming from pairs of points on different $z$-slices (resp. $y$-; $x$-).

\item If $\alpha\ge3$, then $S$ has at least $\frac{1}{64}(\alpha-1)(\alpha-3)^2 n$ \X-edges and \Y-edges total, for any pair of primary colors $\set{\X,\Y} \belongs \set{\R,\G,\B}$.
(If $S$ is slice-increasing in only two of the three coordinate-orientations, the proof still gives something slightly weaker.)

\item If $\alpha\ge3$, then $S$ has at least $\frac{1}{64}(\alpha-1)(\alpha-3)^2 n$ \K-edges.
(This holds even if $S$ is slice-increasing in only two of the three coordinate-orientations.)

\item If $\alpha \ge 100$ and $n \ge M\ge 100n/\alpha$, then there are at least $2^{-11}\alpha^3 M$ \K-edges such that $\Delta y \le M$ while $\Delta x,\Delta z \ge \alpha/8$.
(This holds even if $S$ is only $x$- and $y$- slice-increasing.)
\end{enumerate}
}

We give the proofs below in order.
In all cases the edges we count actually have some additional structure, as the proofs themselves will clarify.
The concrete results do not seem to directly apply to the original problems, since we have not been able to prove any nontrivial \emph{upper} bounds on edge-counts.
But perhaps some of the methods will inspire interested readers.
Curiously, for $\alpha = n^{1/2}$ with $n$ sufficiently large, the first and third bounds are, up to symmetry, sharp up to an absolute constant.

\ex{[Lexicographic bias]
Take for $S$ the image of the map $[m]^3\to \RR^3$ given by $(i,j,k) \mapsto (im+j,im+k,jm+k)$.
Here $\card{S} = m^3$.
It turns out that color classes \R and \K are both dense, with around half of the edges each, while \B and \G are not dense.
Indeed, we check that $(i,j,k) < (i',j',k')$ is colored as follows.
\begin{itemize}
\item \B if $(i',j') = (i,j)$ and $k' > k$, i.e. the first coordinate increase is in position 3.
The \B-edges are all ``trivial'' of type $(0,+,+)$, with $im+j$ fixed.
Total: $m\cdot m\cdot \binom{m}{2} \sim m^4$.

\item \G if $i' = i$ and $j' > j$ and $k' \le k$.
In particular this requires the first coordinate increase to be in position 2.
Total: at most $m\cdot \binom{m}{2}\cdot m^2 \sim m^5$ if we forget the $k'\le k$ condition (in which case we are including some \K's as well, about half-half split).

\item \R if $i' > i$ and $(j',k') \le (j,k)$ lexicographically.
In particular this requires the first coordinate increase to be in position 1.
Total: at most $\binom{m}{2}\cdot m^2\cdot m^2\sim m^6$ if we forget the $(j',k')\le(j,k)$ condition (in which case we are including some \K's as well, about half-half split).

\item \K otherwise.
Most of these come from when $i'>i$ and $(j',k') > (j,k)$ lexicographically, and the other from $i' = i$ and $j'>j$ but $k'>k$.
Total: $\sim m^6$.
\end{itemize}
}

\subsubsection{First bound}

We are counting pairs of points on $z$-slices.
By convexity, this is at least $n\binom{\card{S}/n}{2} = n\binom{\alpha}{2} = \frac12\alpha(\alpha-1)n$, as desired.

\subsubsection{Second bound (using overlapping label-induced corners)}
\label{SUBSUBSEC:L^2-large-implies-edge-lower-bound}

Recall the $L^2$-idea of Section \ref{SUBSEC:L^2-slice-counts-approach}.
Explicitly, $\sum_{i=1}^{n} (\#z=i)^2 = \card{S_z} + \sum \beta(x,y)$, where $\beta(x,y)$ denotes the number of (label-induced) corners overlapped at square $(x,y)\in [n]^2$.
For convenience set $\beta = 0$ for $(x,y)\in S_z$.
Each $(x,y)\notin S_z$ gives rise to $\beta$ ``trivial'' \R-edges (coming from $z$-slices) and $\beta(\beta - 1)/2$ ``nontrivial'' \R-edges.
Furthermore, each $(+,+,?)$-type edge contributes to at most 2 of the corners $(x,y)\notin S_z$.

If $\card{S} = \alpha n$, then by convexity $\sum (\# z = i)^2 \ge \alpha^2 n$, so $\sum \beta \ge \alpha^2 n - \alpha n$.
Now by convexity of the quadratic $u\mapsto \frac12 u(u-1)$ vanishing at $0$, we obtain
\mathd{
\#\text{nontrivial \R-edges from corner overlaps}
\ge 
\frac12 \sum_{(x,y)\notin S_z} \frac12 \beta(\beta - 1)
\ge \frac14 n^2(-1 + (\alpha^2 n - \alpha n)/n^2)^2
}
as long as $(\alpha^2 n - \alpha n)/n^2 \ge 1$.
But $\alpha \ge 3n^{1/2}$ by assumption, and $\alpha\le n$ trivially in general, so $\alpha^2 - \alpha - n \ge \alpha^2 - 2n \ge \alpha^2/2$.
Thus $(\alpha^2 n - \alpha n)/n^2 \ge 1$ and our valid bound simplifies to $(-n + \alpha^2 - \alpha)^2/4 \ge (\alpha^2/2)^2/4 = \alpha^4/16$, as desired.

\subsubsection{Third and fourth bounds}

Consider the $xy$-grid view.
Counting increasing (resp. decreasing) L-shapes gives a lower bound on \K-edges (resp. \B-edges and \G-edges combined).
We will only explicitly write out the proof for increasing L-shapes.
Asymptotically, if $\card{S} = n^{2-\eps}$ and $\alpha = n^c = n^{1-\eps}$, then we get a count on the order of $n^{1+3c} = n^{4-3\eps}$.
In the proof below, one can roughly think of the factor $n$ as choosing a column, the first $n^c$ as choosing a corner for the L, and then $(n^c)^2$ as the two endpoints of the L.

\prop{
\label{PROP:lower-bound-on-K-edges-slice-increasing-sets}

A slice-increasing set $S\belongs [n]^3$ of size $\alpha n$ with $\alpha \ge 3$ has at least $\frac{1}{64}(\alpha-1)(\alpha-3)^2 n$ \K-edges.
}

\pf{
Purge rows with at most $d$ points (where $d$ may be non-integral) to get a new good set $S'$ of size $\beta n \ge \card{S} - dn$, where $\beta \ge \alpha - d$.
If column $x$ of $S'$ has $e = c_x$ points $P_1,\dots,P_e$ (in increasing order), then there must be at least $\sum_{i=1}^{e} (d-1)\cdot \min(i-1,e-i)
\ge (d-1)\cdot \frac14 (e-1)(e-2)$
\K-edges coming from increasing L-shapes whose corner lies on column $S'$: the contribution of $ (d-1)\cdot \min(i-1,e-i)$ comes from corner $P_i$.
On the other hand, each such \K-edge belongs to at most 2 such L-shapes.
We obtain a lower bound
\mathd{
\#\text{\K-edges from L-shapes}
\ge \frac12 \sum_{x=1}^{n} \sum_{i=1}^{c_x} (d-1)\cdot \min(i-1,c_x-i)
\ge \frac12 (d-1)\sum_{x=1}^{n} \frac14 (c_x - 1)(c_x - 2),}
which is at least $\frac12 (d-1)\cdot n\cdot \frac14 (\beta - 1)(\beta - 2)
\ge \frac{n}{8} (d-1)(\alpha - d - 1)(\alpha - d - 2)$ by convexity of the quadratic $u\mapsto \frac14 (u-1)(u-2)$.
Choosing $d = (\alpha-1)/2$ gives a lower bound of $\frac{n}{8}\frac{\alpha-3}{2}\frac{\alpha-1}{2}\frac{\alpha-3}{2}$, as desired.
}

\subsubsection{Fifth bound}

If $\alpha \ge 100$ and $n \ge M\ge 100n/\alpha$, then we can similarly count L-shapes with $\Delta y \le M$ while $\Delta x,\Delta z \ge \alpha/8$.
We will only use that $S$ is $x$- and $y$- slice-increasing.
First purge rows with at most $d = \alpha/2$ points (where $d$ may be non-integral) to get a new set $S'$ of size $\beta n \ge \card{S} - dn$, where $\beta \ge \alpha/2$.

Take a column $x$ of $S'$ with $e = c_x$ points $P_1,\dots,P_e$ (in increasing order).
For each value of $y\in[n]$, let $\nu_y$ (resp. $\eta_y$) denote the number of points $P_k = (x,y',z')$ with $\abs{y'-y} \le M$ and $y'>y$ (resp. $y'<y$).
We estimate $\sum_{y\in[n]: (x,y)\in S_z} \min(\nu_y,\eta_y)$ by finding $u\in [n/M,n/M+1]$ length-$M$ disjoint intervals $I_1,\dots,I_u$ (independent of column $x$) covering $[n]$, containing the $y$-coordinates of exactly $a_1,\dots,a_u$ points among $P_1,\dots,P_e$, respectively, so that $a_1+\dots+a_u = e$.
Then the sum $\sum_{y\in I_j: (x,y)\in S_z} \min(\nu_y,\eta_y)$ is at least $\sum_{i=1}^{a_j} \min(i-1,a_j-i) \ge (a_j-1)(a_j-2)/4$.
Thus our contribution of ``$y$-bounded'' increasing L-shapes with corner on column $x$ is at least
\matha{
\sum_{y\in[n]:(x,y)\in S_z}(d-1-2\cdot \alpha/8)\cdot \min(\nu_y,\eta_y)
&\ge (-1+\alpha/4)\sum_{j=1}^{u} (a_j-1)(a_j-2)/4 \\
&\ge (-1+\alpha/4)u\cdot (-1+e/u)(-2+e/u)/4,
}
where the $1+2\cdot \alpha/8$ correction term excludes possibilities with $\Delta x < \alpha/8$ or $\Delta z < \alpha/8$.
For the final bound, we sum over all columns, and divide by 2 for possible over-counting:
\mathd{
\frac12 \sum_{x=1}^{n} (-1+\alpha/4)u(-1+c_x/u)(-2+c_x/u)/4
\ge \frac{\alpha-4}{32}un\cdot(-1+\beta/u)(-2+\beta/u),
}
i.e. $(\alpha-4)nu^{-1}(\beta-u)(\beta-2u)/32$.
Note that $\beta-2u \ge \alpha/2 - 2n/M - 2 \ge \alpha/4$ since $\alpha\ge100$ and $M \ge 100n/\alpha$ by assumption.
Also, $\alpha-4 \ge \alpha/2$, and $u\le n/M+1 \le 2n/M$ implies $nu^{-1} \ge M/2$.
So $(\alpha-4)nu^{-1}(\beta-u)(\beta-2u)/32 \ge (\alpha/2)(M/2)(\alpha/4)^2/32 = 2^{-11}\alpha^3 M$, as desired.

\end{document}